\pdfoutput=1
\documentclass[11pt]{article}
\usepackage[margin=1in]{geometry}
\usepackage{amsmath,amssymb,amsthm,mathtools,booktabs,array,microtype,placeins,float}
\usepackage{graphicx}
\usepackage{xcolor}
\usepackage[colorlinks=true,linkcolor=black,citecolor=black,urlcolor=black]{hyperref}

\newtheorem{theorem}{Theorem}[section]
\newtheorem{proposition}[theorem]{Proposition}
\newtheorem{lemma}[theorem]{Lemma}
\newtheorem{remark}[theorem]{Remark}
\newcommand{\R}{\mathbb R}
\newcommand{\dd}{\,\mathrm d}
\newcommand{\Id}{\mathrm I}
\newcommand{\erfcx}{\operatorname{erfcx}}
\newcommand{\SL}{\mathcal{S}}
\newcommand{\DL}{\mathcal{D}}

\title{High-Order Regularization for Singular and Nearly Singular Laplace and Stokes Layer Potentials in Two Dimensions}
\author{Han Zhou\thanks{Department of Mathematics, University of Pennsylvania,
Philadelphia, PA, USA (\texttt{hzhou24@sas.upenn.edu}).}}
\date{}
\begin{document}
\maketitle
\begin{abstract}
Layer potentials represent the solutions of many elliptic boundary value
problems, but their kernels have the singularity of the fundamental
solution, so standard quadrature loses accuracy at targets on or near the
boundary. In two dimensions, Gaussian regularization of the logarithmic
kernel gives a smooth kernel involving the exponential integral $E_1$, but
for the single layer its regularization error is only first order in the
smoothing radius $\delta$. We construct modified kernels for two-dimensional
Laplace and Stokes single- and double-layer potentials that achieve
regularization error $O(\delta^p)$, $p=3,5,7$, uniformly for targets on and
near the curve. The resulting method is both high-order accurate and easy to
implement: the integrals are evaluated directly by the periodic trapezoidal
rule on a single fixed set of boundary nodes for all targets. Choosing
$\delta=c h^q$ with $q<1$, where $h$ is the node spacing, keeps the
discretization error smaller than the regularization error as $h\to0$, and
numerical experiments confirm the predicted $O(h^{pq})$ convergence.
Finally, we combine the near-curve evaluation with a fast Poisson solver to
compute the potentials on the entire Cartesian grid.
\end{abstract}

\setcounter{tocdepth}{1}

\section{Introduction}\label{sec:intro}
Boundary integral methods reduce elliptic problems with boundary or interface conditions
to layer potentials on lower-dimensional boundaries. However, the kernels inherit the singularity of the
fundamental solution, and standard quadratures lose accuracy at target points on or close to the boundary.
Evaluation directly on the curve (the singular case) is central to
boundary integral equations \cite{GKM1996,Kress2014}. Evaluation near the curve (the nearly singular
case) arises in interacting drops or vesicles in Stokes flow \cite{OT2015,Pozrikidis2001}, and in
grid-based methods requiring values at Cartesian points adjacent to an interface
\cite{Beale2004,Mayo1984,Mayo1985,MGM1995}. Here we derive high-order regularized kernels for
two-dimensional layer potentials with regularization error $O(\delta^p)$ ($p=3,5,7$), allowing accurate
evaluation by the periodic trapezoidal rule without specialized near-singularity quadratures.
We treat harmonic single- and double-layer potentials, the Stokeslet, stresslet, and the pressure in
Stokes flow. The formulas hold uniformly near and on the boundary.

The nearly singular regime is challenging because, for an analytic periodic integrand, the
convergence rate of the trapezoidal rule is governed by the distance from the real axis to the
nearest singularity in the complex parameter plane \cite{TW2014}. For layer potentials, this distance
is proportional to the physical distance $|b|$ from the target to the boundary $\Gamma$ \cite{Barnett2014}.
Consequently, accuracy deteriorates as $|b|/h \to 0$, where $h$ is the mesh size, even when the geometry
and data are analytic.

To overcome this loss of accuracy, we regularize the integral by replacing the singular kernel
with a smooth one that depends on a smoothing radius $\delta$. In three dimensions, Beale and collaborators
replace $1/r$ by $\operatorname{erf}(r/\delta)/r$ \cite{BJRT2022,BST2024,BYW2016,TB2013,TB2019} and
either extrapolate in $\delta$ \cite{BT2024} or modify the smoothing factor to obtain error $O(\delta^p)$
with $p=3,5,7$ \cite{BT2026}. In two dimensions the fundamental solution is logarithmic, and Gaussian
regularization adds a term to $\log r$ rather than multiplying by a cutoff factor. Smoothed
logarithmic kernels are standard in vortex methods for area integrals \cite{BM1985} and vortex sheet
desingularization \cite{BB2004,Krasny1986}, and were used for boundary potentials by
Beale and Lai \cite{BL2001}. Using the gradient of the regularized kernel, which involves the
exponential integral $E_1$, they treated double layers and tangential derivatives of single layers.
After subtracting the density at the closest point, their double-layer error was $O(\delta^2)$, and the
single-layer error was $O(\delta)$. They eliminated these leading terms with computed corrections
requiring derivatives of the density and parametrization, and corrected the discretization error for
$\delta$ proportional to $h$. For the two-dimensional logarithmic kernel, however, existing kernel
regularizations have not achieved higher-order accuracy uniformly near the boundary without derivatives
of the density. Here we extend the approach of \cite{BT2026} by modifying the $E_1$ kernel to achieve
uniform regularization error $O(\delta^p)$ for $p=3,5,7$.
Higher order permits a larger ratio $\delta/h$, which governs the trapezoidal discretization error.
If $\delta/h$ is fixed as $h\to0$, the discretization error decreases only slowly ($O(h)$ in \cite{BL2001},
corrected to $O(h^2)$). We therefore set $\delta=c h^q$ with $q<1$, ensuring the discretization error
decays faster than the regularization error and yielding convergence $O(h^{pq})$; for example,
$p=7$ and $q=5/7$ give $O(h^5)$. With $\delta$ proportional to $h$, we observe $O(h^p)$ over the range
where discretization error is negligible.

We illustrate the approach with the single-layer potential
\[
 Sf(y)=\int_\Gamma G(|x-y|)f(x)\dd s_x,\qquad G(r)=\frac{1}{2\pi}\log r,
\]
on a smooth closed curve $\Gamma$ with given density $f$. The regularized
version is
\[
 S_\delta^{(1)}f(y)=\int_\Gamma G_\delta^{(1)}(|x-y|)f(x)\dd s_x,\qquad
 G_\delta^{(1)}(r)=G(r)+\frac{1}{4\pi}E_1\bigl(r^2/\delta^2\bigr).
\]
Here $E_1(\xi)=\int_\xi^\infty v^{-1}e^{-v}\dd v$ is the exponential integral.
The function $G_\delta^{(1)}$ is smooth. If $y$ is near $\Gamma$, we can write $y=x_0+bn$, where $x_0$ is the closest
point to $y$ on $\Gamma$, $n$ is the outward normal at $x_0$, and $b$ is the
signed distance. With $\lambda=b/\delta$, we show in
Section~\ref{sec:analysis} that
\[
 S_\delta^{(1)}f(y)=Sf(y)+\frac{1}{4\pi}\Bigl[w_0\delta I_0^{(1)}(\lambda)
 +w_2\delta^3I_2^{(1)}(\lambda)+w_4\delta^5I_4^{(1)}(\lambda)\Bigr]+O(\delta^7),
\]
uniformly for $y$ within a fixed multiple of $\delta$ of $\Gamma$. Here the
$I_n^{(1)}$ are moments that are known explicitly, and the Taylor
coefficients $w_n$ depend on the curve, the density, and $y$, but not on
$\delta$.

To make these moments vanish, we add scaled derivatives of $E_1$, where $E_1^{[k]}$ denotes the
$k$-th derivative,
\begin{equation}\label{eq:ansatzintro}
 G_\delta^{(p)}(r;b)=G(r)+\frac{1}{4\pi}s^{(p)}\!\left(\frac{r^2}{\delta^2};\lambda\right),
 \qquad
 s^{(p)}(\xi;\lambda)=E_1(\xi)+\sum_{k=1}^{m}a_{k}^{(p)}(\lambda)\,\xi^kE_1^{[k]}(\xi),
\end{equation}
with $m=(p-1)/2$ coefficients determined by the conditions that the first
$m$ moments vanish; the resulting regularization error is $O(\delta^p)$. Integration by parts relates
the double-layer and Stokes kernels to the single layer, so that they use the same coefficients
$a_k^{(p)}(\lambda)$, together with an explicit local correction at $x_0$ for one remaining moment in the
Stokes kernels. The pressure is written in terms of Laplace single- and double-layer potentials and uses
the same kernels.

We test the method on Laplace and Stokes problems in several geometries
and verify the predicted orders of accuracy for target points on and near the boundary.
Following Mayo \cite{Mayo1984,Mayo1985}, Beale and Lai \cite{BL2001}, and Beale and Tlupova \cite{BT2026},
we also use the regularized integrals at grid points near the curve to compute solutions on the entire
Cartesian grid with a fast Poisson solver.

The method belongs to a family of approaches that analyze the local error
of a smoothed or trapezoidal approximation and remove its leading terms.
For evaluation on a curve, corrected trapezoidal rules for logarithmic and
other singularities \cite{Alpert1999,KR1997,SI1988,WM2021} and the product
quadrature of Kress \cite{Kress2014} achieve high order. Near the curve,
related methods include the local correction of \cite{BL2001}, the
asymptotic approximations of \cite{CKK2018}, and the corrected trapezoidal
rules for near-singular integrals of \cite{Nitsche2022}. In three
dimensions, related constructions include the extrapolated regularization
of \cite{BT2024}, its fast-summation variant \cite{SieborTlupova2025}, the
high-order smoothing functions of \cite{BT2026} and \cite{FPT2026}, and
corrected trapezoidal rules for near-singular integrals on surfaces
\cite{NWX2025}; Richardson extrapolation has also been applied to
regularized Stokeslets \cite{GS2021}. Conditions under which
discretization corrections can be neglected are analyzed in
\cite{Beale2020}.

Other high-order approaches treat the local singularity directly. In two
dimensions, Helsing and Ojala \cite{HO2008} evaluate layer potentials close
to the curve by panel-based product integration, an approach extended to
drops in Stokes flow in \cite{OT2015}, and Barnett and coauthors construct
quadratures that remain accurate close to the curve for Laplace, Helmholtz,
and Stokes potentials \cite{Barnett2014,BWV2015}. Helsing's higher-order
singularity subtraction \cite{Helsing2013} expands singular and nearly singular kernels
locally and integrates the resulting terms by high-order product integration on
parametrized surface patches; quadrature by
expansion \cite{KBG2013} evaluates the potential from an expansion about a
point away from the boundary; harmonic density interpolation \cite{PFT2019}
replaces the density locally by a function whose layer potential is known
analytically, and \cite{FPB2021} extends density interpolation to other
elliptic operators. Values near a surface can also be interpolated from
points farther away \cite{YBZ2006}, and quadrature by fundamental solutions
\cite{SB2022} provides a kernel-independent alternative for simple objects.
The method of regularized Stokeslets \cite{Cortez2001,CFM2005} also
replaces the singular kernel by a smooth one, but it uses a fixed smoothed
source and does not cancel moments of the error; its regularization error
can be reduced for rigid objects in three dimensions \cite{NC2014}.

A practical advantage of the method is its simplicity. Once the regularized
kernels are defined, potentials are evaluated directly with the periodic trapezoidal rule on a fixed set of
boundary nodes for all targets, provided $\delta/h$ is sufficiently large. The regularized kernel differs
appreciably from the singular kernel only within a neighborhood of width $O(\delta)$ around the target.

Section~\ref{sec:e1ansatz} presents the regularized Laplace kernels and uniform error estimates, and
Section~\ref{sec:stokes} treats the Stokes kernels and pressure. Section~\ref{sec:analysis} derives the
asymptotic expansions, proves the error estimates, and bounds the discretization error. Numerical tests
are given in Sections~\ref{sec:num} and~\ref{sec:num-allgrid}, conclusions in
Section~\ref{sec:conclusions}, and Appendix~\ref{app:implementation} gives stable formulas for
evaluating the kernels and coefficients without cancellation.

\section{Regularized kernels for the Laplace equation}\label{sec:e1ansatz}
\subsection{The single layer and the \texorpdfstring{$E_1$}{E1} kernel}\label{sec:gauss2d}
Let $\Gamma\subset\R^2$ be a smooth, simple, closed curve, oriented
counterclockwise, with outward unit normal $n$, unit tangent $\tau$,
arclength $s$, and signed curvature $\kappa$, defined by $\tau_s=-\kappa n$ and
$n_s=\kappa\tau$, so that $\kappa>0$ on a convex curve. Write $\Omega$ for the interior domain. For a target $y$ in a fixed,
sufficiently narrow tubular neighborhood of $\Gamma$ where the closest boundary point is unique,
let $x_0=x_0(y)$ be that closest point and write
\begin{equation}\label{eq:offset}
 y=x_0+b\,n_0,\qquad n_0=n(x_0),\qquad \lambda=\frac{b}{\delta}.
\end{equation}
Here $\delta$ is the smoothing radius, and $\lambda$ is the scaled signed
distance from $\Gamma$ to the target, with $b>0$ outside $\Omega$.

The source point $x$ lies on $\Gamma$, and $y$ is the target point where
the potential is evaluated. We set $\boldsymbol r=x-y=(r_1,r_2)$ and
$r=|\boldsymbol r|$. We regularize the logarithmic fundamental solution
\begin{equation}\label{eq:logkernel}
 G(r)=\frac{\log r}{2\pi}.
\end{equation}
The single-layer potential with smooth density $f$ on $\Gamma$ is
\begin{equation}\label{eq:singlelayer}
 (Sf)(y)=\int_\Gamma G(|x-y|)f(x)\dd s_x .
\end{equation}
Because the double-layer and Stokes kernels are formed from derivatives of
\eqref{eq:logkernel}, we obtain their regularized versions by differentiating the regularized version of $G$.

Let
\begin{equation}\label{eq:blob}
 \varphi_\delta(r)=\frac{1}{\pi\delta^2}e^{-r^2/\delta^2},
 \qquad \int_{\R^2}\varphi_\delta\dd A=1,
\end{equation}
and set $G_\delta=G*\varphi_\delta$. Since $\Delta G$ is the Dirac delta at the
origin, we have $\Delta G_\delta=\varphi_\delta$; thus $G_\delta$ is the
potential of a radially symmetric Gaussian source of unit mass and width $\delta$.

By radial symmetry and the divergence theorem,
$\nabla G_\delta$ equals $\nabla G$ times the mass of $\varphi_\delta$ inside the
disk of radius $r$, namely
$2\pi\int_0^r v\varphi_\delta(v)\dd v=1-e^{-r^2/\delta^2}$. Thus, as in
\cite[eq.~(1.6)]{BL2001},
\begin{equation}\label{eq:mollifiedgrad}
 \nabla G_\delta=\bigl(1-e^{-r^2/\delta^2}\bigr)\nabla G,
 \qquad
 \partial_rG_\delta=\frac{1}{2\pi r}\bigl(1-e^{-r^2/\delta^2}\bigr).
\end{equation}
Integrating \eqref{eq:mollifiedgrad} in $r$ and choosing the constant of integration so
that $G_\delta-G\to0$ as $r\to\infty$ gives the regularized
kernel
\begin{equation}\label{eq:Greg0}
 G_\delta^{(1)}(r)=G(r)+\frac{1}{4\pi}E_1(\xi),\qquad
 \xi=\rho^2,\qquad \rho=\frac r\delta,\qquad
 E_1(\xi)=\int_\xi^\infty\frac{e^{-v}}{v}\dd v . 
\end{equation}
In three dimensions,
the analogous calculation gives $\operatorname{erf}(\rho)/r$, which differs from
$1/r$ by a term with Gaussian decay; in two dimensions, it gives $\log r/(2\pi)$
plus a multiple of $E_1$. The function $E_1$ is available in standard
mathematical libraries; for small $r/\delta$ we evaluate $G_\delta^{(1)}$ from its
power series (Appendix~\ref{app:implementation}).

For $y$ near $\Gamma$, we can expand the regularization error in powers
of $\delta$. In Section~\ref{sec:baserror} we show that
\begin{equation}\label{eq:baseexp}
 \begin{aligned}
 S_\delta^{(1)}f-Sf
 &\sim\frac{1}{4\pi}\Bigl[w_0\delta I_0^{(1)}(\lambda)
 +w_2\delta^3I_2^{(1)}(\lambda)+w_4\delta^5I_4^{(1)}(\lambda)+\cdots\Bigr],
 \\
 I_n^{(1)}(\lambda)&=\int_\R t^nE_1(t^2+\lambda^2)\dd t,
 \qquad n=0,2,4,\ldots .
 \end{aligned}
\end{equation}
Here $\eta$ is a coordinate along $\Gamma$ with $|x-y|^2=\eta^2+b^2$, of the
same sign as the arclength from $x_0$, and $w_n(x_0,b)$ are the even Taylor
coefficients in $\eta$ of the density times $\dd s/\dd\eta$. They depend on
the curve, the density, and the target distance $b$, but not on the
regularized kernel; the kernel enters only through the moments
$I_n^{(1)}(\lambda)$. The power $\delta^{n+1}$ comes from $\eta^n=\delta^nt^n$
and $\dd\eta=\delta\dd t$. Since $E_1$ is positive and decreasing, the first moment
$I_0^{(1)}(\lambda)$ is positive and decreases with $|\lambda|$. Hence, whenever $w_0\neq0$, the $E_1$ kernel is only first
order in $\delta$, and this error cannot be reduced by taking more quadrature points.

\subsection{Higher-order kernels}\label{sec:highorder}
To get higher order, we modify the $E_1$ term so that the leading moments in
\eqref{eq:baseexp} are zero. Since successive terms of \eqref{eq:baseexp}
differ by a factor $\delta^2$, cancelling each moment increases the
order by two and requires one free coefficient in the
kernel.
The kernel \eqref{eq:Greg0} is smooth because the logarithmic singularities
of $G$ and $E_1$ cancel at $r=0$. A term added to $E_1$ must keep this
cancellation and decay exponentially, so that the regularized kernel differs
from $G$ only near the source; a multiple of $E_1$ would spoil the
cancellation. The scaled derivatives $\xi^kE_1^{[k]}$, $k\ge1$, where
$E_1^{[k]}=\dd^kE_1/\dd\xi^k$, have both properties: each is a polynomial
times $e^{-\xi}$ (Section~\ref{sec:coeff}). They are easy to differentiate for
the double-layer and Stokes kernels, and their moments are known in closed
form.

Because the moments depend on $\lambda$, the moment conditions must hold for each
$\lambda$ in the band of interest, which cannot be achieved
with constant coefficients \cite[Sect.~2]{BL2001}. We
therefore allow the coefficients to depend on $\lambda$ and, for order
$p=2m+1$, define
\begin{equation}\label{eq:F}
 s^{(p)}(\xi;\lambda)=E_1(\xi)+\sum_{k=1}^{m}a_{k}^{(p)}(\lambda)\,\xi^kE_1^{[k]}(\xi),
\end{equation}
with $a_{k}^{(p)}(\lambda)=0$ for $k>m$. The regularized
single-layer kernel is
\begin{equation}\label{eq:Greg}
 G_\delta^{(p)}(r;b)=G(r)+\frac{1}{4\pi}s^{(p)}(\xi;\lambda),\qquad
 \xi=\frac{r^2}{\delta^2}, 
\end{equation}
where $p=1$ gives the $E_1$ kernel of Section~\ref{sec:gauss2d}, and $p=3,5,7$
give third-, fifth-, and seventh-order kernels. We call $s^{(p)}$ the regularizing
function. Although $s^{(p)}$ has a logarithmic singularity at $\xi=0$, its sum with $G$
in \eqref{eq:Greg} is smooth. For a fixed target $y$, $\lambda$ is treated as a constant
when we differentiate the kernel with respect to $x$. Note that the coefficients depend only on
$\lambda$, and not on the density, the curve, or the mesh; this
dependence on $\lambda$ is what distinguishes our kernel from a
convolution with a fixed blob. The regularized
single-layer potential is
\begin{equation}\label{eq:Sreg}
 S_\delta^{(p)}f(y)=\int_\Gamma G_\delta^{(p)}(|x-y|;b)\,f(x)\dd s_x,
\end{equation}
with $b$ and $\lambda=b/\delta$ defined in \eqref{eq:offset}.

The logarithmic singularity cancels for every $p$. Since
$E_1(\xi)=-\gamma-\log\xi+\xi-\cdots$, where $\gamma$ is Euler's constant, we
have $\xi^kE_1^{[k]}(\xi)\to(-1)^k(k-1)!$ as $\xi\downarrow0$ for $k\ge1$, and
the value of the regularized kernel at $r=0$ is finite:
\begin{equation}\label{eq:coincident}
 G_\delta^{(p)}(0;0)=\frac{1}{4\pi}
 \Bigl[2\log\delta-\gamma+\sum_{k=1}^{m}(-1)^k(k-1)!\,a_{k}^{(p)}(0)\Bigr].
\end{equation}
Since $r\ge|b|$, we have $r=0$ only when the target lies on the curve
($b=\lambda=0$) and $x=y$; using \eqref{eq:coincident} directly
avoids subtracting two large numbers (with $a_k^{(p)}(0)$ from
\eqref{eq:c3}--\eqref{eq:c7}, the sum in \eqref{eq:coincident} equals $-2$,
$-\tfrac83$, and $-\tfrac{46}{15}$ for $p=3,5,7$).

The same expansion for the kernel \eqref{eq:Greg} gives
\begin{equation}\label{eq:SLerror}
 S_\delta^{(p)}f-Sf\sim\frac{1}{4\pi}
 \sum_{\substack{n\ge0\\ n\ {\rm even}}}
 w_n(x_0,b)\,\delta^{n+1}I_n^{(p)}(\lambda),
 \qquad
 I_n^{(p)}(\lambda)=\int_\R t^n
 s^{(p)}(t^2+\lambda^2;\lambda)\dd t .
\end{equation}
Here only the moments change, while the Taylor coefficients $w_n$ are the same as in
\eqref{eq:baseexp}. For the first $m$ terms of
\eqref{eq:SLerror} to vanish for arbitrary curves and densities, the first $m$ even
moments must be zero:
\begin{equation}\label{eq:M}
 I_n^{(p)}(\lambda)=0,\qquad n=0,2,\dots,p-3.
\end{equation}
This gives $m$ linear equations for the $m$ coefficients $a_k^{(p)}(\lambda)$ in \eqref{eq:F}. When
\eqref{eq:M} holds, the leading nonzero term in \eqref{eq:SLerror} is $n=p-1$,
which is $O(\delta^{p})$ because $I_{p-1}^{(p)}$ and $w_{p-1}$ are uniformly bounded
for $|b|\le\Lambda\delta$ (with $\Lambda$ fixed) when the curve and
density are smooth. Thus the kernel \eqref{eq:Greg} has $O(\delta^p)$ regularization error.

For $p=3,5,7$, the linear system \eqref{eq:M} can be solved explicitly
(Section~\ref{sec:coeff}) in terms of the integrals
\begin{equation}\label{eq:Hdef}
 \mathcal{H}_n(\lambda)=\int_1^\infty e^{-\lambda^2(v-1)}v^{-(n+3)/2}\dd v,
 \qquad \mathcal{H}_n(0)=\frac{2}{n+1},\qquad n=0,2,4,\ldots,
\end{equation}
which are positive, bounded, even functions of $\lambda$ that decrease with
$|\lambda|$. Appendix~\ref{app:coeff} gives $\mathcal{H}_0$ in closed form and a
recurrence for the other $\mathcal{H}_n$, and describes how to evaluate them and
the coefficients accurately for all $\lambda$. In terms of $\mathcal{H}_n$, we have
$I_n^{(1)}(\lambda)=\Gamma((n+1)/2)e^{-\lambda^2}\mathcal{H}_n(\lambda)$ for
even $n\ge0$, and the coefficients are
\begin{align}
 p=3:\quad &a_{1}^{(3)}=\mathcal{H}_0,\label{eq:c3}\\[2pt]
 p=5:\quad &a_{2}^{(5)}=\mathcal{H}_0-\mathcal{H}_2,\qquad
 a_{1}^{(5)}=\mathcal{H}_0+\bigl(\lambda^2+\tfrac32\bigr)a_{2}^{(5)},\label{eq:c5}\\[2pt]
 p=7:\quad &a_{3}^{(7)}=\tfrac12(\mathcal{H}_0-2\mathcal{H}_2+\mathcal{H}_4),\notag\\
 &a_{2}^{(7)}=\mathcal{H}_0-\mathcal{H}_2+\bigl(2\lambda^2+5\bigr)a_{3}^{(7)},\notag\\
 &a_{1}^{(7)}=\mathcal{H}_0+\bigl(\lambda^2+\tfrac32\bigr)a_{2}^{(7)}
 -\bigl(\lambda^4+3\lambda^2+\tfrac{15}{4}\bigr)a_{3}^{(7)}.\label{eq:c7}
\end{align}
On the curve ($\lambda=0$), these take the simple values $a_{1}^{(3)}=2$,
$(a_{1}^{(5)},a_{2}^{(5)})=(4,\tfrac43)$, and
$(a_{1}^{(7)},a_{2}^{(7)},a_{3}^{(7)})=(6,4,\tfrac{8}{15})$. Just as in
\cite{BT2026}, the formulas are in triangular form: setting $a_3^{(7)}=0$ in the $p=7$ formulas gives the coefficients
for $p=5$, and setting also $a_2^{(5)}=0$ gives $a_1^{(3)}$.

With $p=3$, the kernel cancels the leading error term that Beale and Lai
compute and subtract \cite{BL2001}. Their correction requires derivatives of
the density; the modified kernel does not, so it is simpler to use, and
$p=5$ or $7$ cancels further terms in the same way.

\subsection{Radial factors}\label{sec:props}
To obtain the double-layer and Stokes kernels, we write the radial
derivative of the regularized kernel as $\partial_r G$ multiplied by a radial
factor:
\begin{equation}\label{eq:defectrelation}
 \partial_r\Bigl[G(r)+\frac{1}{4\pi}s^{(p)}(\xi;\lambda)\Bigr]
 =\frac{1}{2\pi r}\bigl(1+\xi\partial_\xi s^{(p)}(\xi;\lambda)\bigr),
 \qquad \xi=\frac{r^2}{\delta^2},
\end{equation}
so the factor is $1-Q^{(p)}$ with $Q^{(p)}=-\xi\partial_\xi s^{(p)}$. For $p=1$, we have
$Q^{(1)}=e^{-\xi}$ by \eqref{eq:mollifiedgrad}; adding the scaled derivatives
of $E_1$ multiplies this Gaussian by a polynomial of degree
$m$. Specifically, $-\xi\partial_\xi E_1=e^{-\xi}$ and
$-\xi\partial_\xi[\xi^kE_1^{[k]}]=(-1)^k\xi^ke^{-\xi}$ for $k\ge1$
(Section~\ref{sec:coeff}). Differentiating \eqref{eq:F} term by
term gives
\begin{equation}\label{eq:PQT}
 \begin{aligned}
 Q^{(p)}(\xi;\lambda)&:=-\xi\partial_\xi s^{(p)}(\xi;\lambda)
 =e^{-\xi}P^{(p)}(\xi;\lambda),\\
 P^{(p)}(\xi;\lambda)&\mathrel{:=}1+\sum_{k=1}^{m}(-1)^ka_{k}^{(p)}(\lambda)\xi^k.
 \end{aligned}
\end{equation}
We then define the radial factors
\begin{equation}\label{eq:dfactors}
 \begin{aligned}
 d^{(p)}(\xi;\lambda)&:=1-Q^{(p)}(\xi;\lambda)
 =1+\xi\partial_\xi s^{(p)}(\xi;\lambda),\\
 \tilde d^{(p)}(\xi;\lambda)&:=d^{(p)}(\xi;\lambda)
 -\xi\partial_\xi d^{(p)}(\xi;\lambda)
 =1-\xi^2\partial_\xi^2s^{(p)}(\xi;\lambda).
 \end{aligned}
\end{equation}
Thus the single-layer kernel is regularized by adding $s^{(p)}/(4\pi)$, while the
double-layer kernel is $\nabla G$ times $d^{(p)}$. The second factor
$\tilde d^{(p)}$ is used only in the stresslet kernel.
This construction is the
two-dimensional counterpart of \cite[eq.~(13)]{BT2026}, where the three-dimensional smoothing factor
$\operatorname{erf}\rho$ is modified by adding multiples
of $\rho(\operatorname{erf}\rho)'$, $\rho^2(\operatorname{erf}\rho)''$, and
higher scaled derivatives.

Since the coefficients \eqref{eq:c3}--\eqref{eq:c7} are polynomials in
$\lambda^2$ and the bounded functions $\mathcal{H}_n$, they remain bounded for
$|\lambda|\le\Lambda$ for any fixed $\Lambda>0$. Consequently,
$Q^{(p)}=e^{-\xi}P^{(p)}$ decays like a Gaussian uniformly for
$|\lambda|\le\Lambda$, and $G_\delta^{(p)}-G$ is exponentially small when $\rho\gg1$. The
local expansion \eqref{eq:SLerror} therefore holds for all the
regularized kernels $G_\delta^{(p)}$.

\subsection{The double-layer potential}\label{sec:dl}
The double-layer kernel is the normal derivative of $G$
with respect to the source point $x$:
\begin{equation}\label{eq:Dkernel}
 K(x,y)=\partial_{n(x)}G(r)=\frac{n(x)\cdot\boldsymbol r}{2\pi r^{2}},
 \qquad \boldsymbol r=x-y,
\end{equation}
and the double-layer potential with smooth density $g$ is
\begin{equation}\label{eq:Dpot}
 (Dg)(y)=\int_\Gamma K(x,y)g(x)\dd s_x .
\end{equation}
When the target $y$ lies on $\Gamma$, $Dg$ stands for the direct value $\mathcal{K}g$ of the
integral, which equals the average of the interior and exterior limits $D^-g$ and $D^+g$; we
use the same convention for the other potentials that jump across $\Gamma$. The one-sided limits
are given by the jump relations \cite{Kress2014}
\begin{equation}\label{eq:jump}
 D^-g=\mathcal{K}g+\tfrac12g,\qquad D^+g=\mathcal{K}g-\tfrac12g.
\end{equation}
Our singularity subtraction uses Gauss's identity
\begin{equation}\label{eq:chi}
 \int_\Gamma K(x,y)\dd s_x=\chi(y),\qquad
 \chi(y)=\begin{cases}1,&y\in\Omega,\\ 0,&y\notin\overline\Omega,\\
 \tfrac12,&y\in\Gamma.\end{cases}
\end{equation}
Differentiating \eqref{eq:Greg} with respect to $x$ while keeping
$y$, $b$, $\delta$, and thus $\lambda$ and $a_k^{(p)}(\lambda)$ fixed gives the regularized double-layer kernel
\begin{equation}\label{eq:Kreg}
 K_\delta^{(p)}(x,y;b)=\partial_{n(x)}G_\delta^{(p)}
 =\frac{n(x)\cdot(x-y)}{2\pi|x-y|^{2}}\,
 d^{(p)}\!\left(\frac{|x-y|^2}{\delta^2};\lambda\right).
\end{equation}
Because $d^{(p)}(0;\lambda)=0$ and $d^{(p)}(\xi;\lambda)/\xi$ is analytic at
$\xi=0$ with value $1+a_{1}^{(p)}$, the kernel \eqref{eq:Kreg} is a smooth function of
$x$ and equals zero when $x=y$. We define the regularized double-layer potential with subtraction by
\begin{equation}\label{eq:Dreg}
 D_\delta^{(p)}g(y)=\int_\Gamma K_\delta^{(p)}(x,y;b)\,[g(x)-g(x_0)]\dd s_x
 +\chi(y)\,g(x_0). 
\end{equation}
For $y\in\Gamma$, \eqref{eq:Dreg} with $\chi(y)=\tfrac12$ approximates the direct value
$\mathcal{K}g$, while setting $\chi=1$ or $0$ gives the interior or exterior limit in \eqref{eq:jump}.
Because of \eqref{eq:chi}, \eqref{eq:Dreg} is exact when $g$ is constant: the
integral is zero and $\chi g$ reproduces the exact potential. Subtracting $g(x_0)$ also makes the
integrand vanish at the closest point $x_0$, removing the constant term in the even part of
$n(x)\cdot(x-y)\,[g(x)-g(x_0)]$. The same subtraction is used in three
dimensions in \cite[Sect.~2]{BT2026} and \cite{TB2019}.

The error $D_\delta^{(p)}g-Dg$ has a local expansion of the same form as
\eqref{eq:SLerror}, with the moments
\begin{equation}\label{eq:Jdef}
 J_n^{(p)}(\lambda)=\int_\R\frac{t^{n+2}}{t^2+\lambda^2}
 Q^{(p)}(t^2+\lambda^2;\lambda)\dd t,\qquad n=0,2,4,\ldots,
\end{equation}
in place of $I_n^{(p)}$. Section~\ref{sec:dlanalysis} derives this expansion
and shows that $J_n^{(p)}=\frac{n+1}{2}I_n^{(p)}$. Therefore, the same
coefficients $a_k^{(p)}(\lambda)$ give $O(\delta^p)$ regularization error for the double layer as well.
Without subtracting $g(x_0)$, the error would not go to zero uniformly as $y$ approaches the
curve, because $K_\delta^{(p)}$ does not satisfy Gauss's identity \eqref{eq:chi}.

\subsection{Error estimate}\label{sec:lapbound}
For $p=3,5,7$, using the coefficients \eqref{eq:c3}--\eqref{eq:c7}, the regularized potentials have
error $O(\delta^p)$ uniformly for targets within a fixed
multiple of $\delta$ of $\Gamma$. The proof is given in
Section~\ref{sec:proofs}.
\begin{proposition}[Uniform accuracy near the curve]\label{prop:lap}
Let $\Gamma$ and the densities $f$ and $g$ be smooth, and let $\Lambda>0$
and $\delta_0>0$ satisfy $\Lambda\delta_0\|\kappa\|_\infty\le\frac12$, with
$\Lambda\delta_0$ small enough that every point within distance $\Lambda\delta_0$
of $\Gamma$ has a unique closest point on $\Gamma$. Let $p=3$, $5$, or $7$,
and let $S_\delta^{(p)}f$ and $D_\delta^{(p)}g$ be given by \eqref{eq:Sreg} and
\eqref{eq:Dreg}, with the kernels \eqref{eq:Greg} and \eqref{eq:Kreg} and the
coefficients \eqref{eq:c3}--\eqref{eq:c7} at $\lambda=b/\delta$.
Then for $0<\delta\le\delta_0$,
\begin{equation}\label{eq:lapbound}
 \sup_{x_0\in\Gamma,\ |b|\le\Lambda\delta}
 \Bigl(|S_\delta^{(p)}f-Sf|+|D_\delta^{(p)}g-Dg|\Bigr)\le C\delta^{p}.
\end{equation}
The constant $C$
depends on finitely many derivatives of $\Gamma$ and the densities, and on
$p$ and $\Lambda$, but not on $\delta$, $b$, or $x_0$.
\end{proposition}

For targets farther from the curve ($|\lambda|$ large), the modification to
the kernel is exponentially small in $|b|/\delta$ and can be neglected. Where
the curvature is large, the condition $\Lambda\delta_0\|\kappa\|_\infty\le\frac12$
requires $\delta$ (and thus $h$) to be small so that $\Lambda$ can be chosen
large enough; a small $h$ is needed in any case to resolve a region of large
curvature.

\begin{remark}\label{rem:oncurve}
Unlike the single-layer kernel, which has a logarithmic singularity at $x=y$,
the double-layer kernel $K(x,y)$ of \eqref{eq:Dkernel} for $y\in\Gamma$ is a
smooth function of $x\in\Gamma$, with limit $\kappa(y)/(4\pi)$ as $x\to y$
\cite{Kress2014}. With this diagonal value, the periodic trapezoidal rule for
$Dg(y)$ converges faster than any power of $h$ \cite{Kress2014,TW2014}. Hence
no regularization is needed on the curve, though \eqref{eq:Dreg} with
$\lambda=0$ remains valid there. This differs from three dimensions, where the
double-layer kernel is $O(1/r)$ on the surface and is regularized there as well
\cite{Beale2004,TB2019}.
\end{remark}

\section{Regularized kernels for Stokes flow}\label{sec:stokes}
We now extend the regularization to two-dimensional Stokes flow, regularizing
the Stokeslet and stresslet integrals and expressing the pressure in terms of harmonic
layer potentials, using the regularized kernels and coefficients of Section~\ref{sec:e1ansatz}.

\subsection{The Stokeslet and stresslet integrals}\label{sec:stokernels}
For creeping flow with viscosity $\mu>0$, velocity $u$ and pressure $p$ satisfy
$\mu\Delta u=\nabla p$ and $\nabla\cdot u=0$, with hydrodynamic stress
$\Sigma=-p\Id+\mu(\nabla u+\nabla u^T)$ \cite{Pozrikidis1992}. From here on, $p$
has two meanings, distinguished by use: the pressure is a field, as in
$\nabla p$, $p^\pm$, and $p_\delta$, and the order is a number, as in the
superscript $(p)$, $\delta^p$, and $p=3,5,7$. The single-layer (Stokeslet)
kernel $S_{ij}$ and double-layer (stresslet) kernel $T_{ijk}$ are
\begin{equation}\label{eq:stokeskernels}
 S_{ij}(\boldsymbol r)=\frac{1}{4\pi\mu}\Bigl(-\delta_{ij}\log r
 +\frac{r_ir_j}{r^2}\Bigr),
 \qquad
 T_{ijk}(\boldsymbol r)=\frac{1}{\pi}\frac{r_ir_jr_k}{r^4},
 \qquad \boldsymbol r=x-y,
\end{equation}
with
\begin{equation}\label{eq:stokeslayers}
 (\SL f)_i(y)=\int_\Gamma S_{ij}(x-y)f_j(x)\dd s_x,\qquad
 (\DL q)_i(y)=\int_\Gamma T_{ijk}(x-y)q_j(x)n_k(x)\dd s_x ,
\end{equation}
where $f$ and $q$ are prescribed smooth vector densities and repeated
indices are summed. We regularize both kernels by expressing them through
$G$ and its radial derivatives and replacing these by their regularized
versions from Section~\ref{sec:e1ansatz}. This gives
\begin{equation}\label{eq:Udelta}
 S_{ij,\delta}^{(p)}
 =-\frac{\delta_{ij}}{2\mu}G_\delta^{(p)}
 +\frac{1}{4\pi\mu}\frac{r_ir_j}{\delta^2}A^{(p)}(\xi;\lambda),
 \qquad
 A^{(p)}(\xi;\lambda)=\frac{d^{(p)}(\xi;\lambda)}{\xi},
\end{equation}
\begin{equation}\label{eq:Hdelta}
 T_{ijk,\delta}^{(p)}
 =\frac{1}{\pi}\frac{r_ir_jr_k}{\delta^4}
 B^{(p)}(\xi;\lambda),
 \qquad
 B^{(p)}(\xi;\lambda)=\frac{\tilde d^{(p)}(\xi;\lambda)}{\xi^2}
 =-\partial_\xi A^{(p)}(\xi;\lambda),
\end{equation}
where all derivatives in $\xi$ are taken at fixed $\lambda$. Both kernels have
finite limits as $r\to0$. Since
$A^{(p)}\to1+a_{1}^{(p)}$ and
$B^{(p)}\to\tfrac12+a_{1}^{(p)}+a_{2}^{(p)}$ as $\xi\downarrow0$, we have
$T_{ijk,\delta}^{(p)}=0$ and
$S_{ij,\delta}^{(p)}=-\delta_{ij}G_\delta^{(p)}(0;0)/(2\mu)$ at $r=0$, with
$G_\delta^{(p)}(0;0)$ as in \eqref{eq:coincident}.
We evaluate the kernels through $A^{(p)}$ and $B^{(p)}$, so that $r^{-2}$ and $r^{-4}$
are never formed; see Appendix~\ref{app:implementation}.

The Stokeslet has an additional moment that the conditions \eqref{eq:M} do
not cancel; it is removed by a local correction
(Section~\ref{sec:stokeslet}). The corrected Stokes single layer is
\begin{equation}\label{eq:SLcorr}
 (\SL_{\delta,{\rm corr}}^{(p)}f)_i=\int_\Gamma S_{ij,\delta}^{(p)}f_j\dd s
 +\frac{|b|\,\mathcal{N}^{(p)}(\lambda)}{4\pi\mu\,\alpha}\,(n_0)_i\,(n_0\cdot f_0),
\end{equation}
where $x_0$, $n_0$, and $b$ are as in \eqref{eq:offset}, $f_0=f(x_0)$,
$\alpha=\sqrt{1+b\kappa_0}$, $\kappa_0$ is the curvature at $x_0$, and
\begin{equation}\label{eq:Lclosed}
 \mathcal{N}^{(p)}(\lambda)=\pi\operatorname{erfc}|\lambda|
 +\sqrt\pi\,|\lambda| e^{-\lambda^2}
 \Bigl[-a_{1}^{(p)}+a_{2}^{(p)}\bigl(\lambda^2+\tfrac12\bigr)
 -a_{3}^{(p)}\bigl(\lambda^4+\lambda^2+\tfrac34\bigr)\Bigr],
\end{equation}
with the coefficients $a_k^{(p)}$ of \eqref{eq:c3}--\eqref{eq:c7} and $a_{k}^{(p)}=0$ for $k>m$.
The function $\mathcal{N}^{(p)}$ is bounded and equals $\pi$ at
$\lambda=0$. In three dimensions \cite[eq.~(41)]{BT2026}, the same moment is
treated differently: $(n_0\cdot f_0)\,n$ is subtracted from the density, which
does not change the Stokeslet integral since $\int_\Gamma S_{ij}n_j\dd s=0$.
Here the density is not modified; because $\mathcal{N}^{(p)}$ has the closed
form \eqref{eq:Lclosed}, the error term is computed and added directly.

With the sign convention of \eqref{eq:stokeskernels}, the analogue of
Gauss's identity \eqref{eq:chi} is
\[
 \int_\Gamma T_{ijk}(x-y)n_k(x)\dd s_x=\chi(y)\delta_{ij}.
\]
We use this identity for singularity subtraction. The moments of the
stresslet expansion are multiples of the $J_n^{(p)}$ of \eqref{eq:Jdef}, except
the lowest one, which is $\tfrac12$ times the Stokeslet moment left
uncancelled (Section~\ref{sec:stresslet}). Adding the corresponding local
correction gives the corrected double layer
\begin{equation}\label{eq:DLcorr}
 (\DL_{\delta,{\rm corr}}^{(p)}q)_i
 =\int_\Gamma T_{ijk,\delta}^{(p)}(q_j-q_{0,j})n_k\dd s+\chi q_{0,i}
 +\frac{|b|\,\mathcal{N}^{(p)}(\lambda)}{2\pi}\,\mathcal{V}_i,
\end{equation}
where $q_0=q(x_0)$ and
\begin{equation}\label{eq:C2}
 \mathcal{V}=\alpha^{-3}\bigl\{n_0(\tau_0\cdot q_s)+\tau_0(n_0\cdot q_s)\bigr\}
 -\frac{b}{2\alpha^3}n_0(n_0\cdot q_{ss})
 +\frac{b^2\kappa_s}{2\alpha^5}n_0(n_0\cdot q_s).
\end{equation}
Here all quantities are evaluated at $x_0$, $\kappa_s$ is the arclength
derivative of the curvature, and $q_s$, $q_{ss}$ denote componentwise first
and second arclength derivatives of the density.

Both corrections vanish at $b=0$, so no correction is needed for targets on
the curve. In the computations, $a_{k}^{(p)}$ is obtained from
\eqref{eq:c3}--\eqref{eq:c7}, $\mathcal{N}^{(p)}$ from \eqref{eq:Lclosed},
and the vector $\mathcal{V}$ directly from \eqref{eq:C2}.
Density values and derivatives at $x_0$ are obtained from the trigonometric
interpolants of the nodal values; arclength derivatives are computed at the
nodes by Fourier differentiation, as $|x'(t)|^{-1}\partial_t$, and
$q_{ss}$ by differentiating the computed nodal values of $q_s$.

\begin{remark}\label{rem:divfree}
The regularized integrals approximate the layer potentials for evaluation, rather than
representing exact Stokes flows with smoothed forces as in the method of regularized Stokeslets
\cite{Cortez2001}. For fixed $\delta>0$, the regularized velocity kernel is not divergence-free in general
(holding $\lambda$ fixed gives $\partial_iS_{ij,\delta}^{(p)}=r_j\partial_\xi d^{(p)}(\xi;\lambda)/(2\pi\mu\delta^2)$),
and the local corrections do not restore incompressibility at finite $\delta$; the divergence-free
condition is recovered to $O(\delta^p)$ as $\delta\to 0$. The pressure is computed
separately from the harmonic potentials $S$ and $D$.
\end{remark}

\subsection{The pressure}\label{sec:pressure}

We express the pressure for the Stokes layer potentials in terms of the harmonic potentials
$S$ and $D$ of \eqref{eq:singlelayer} and \eqref{eq:Dpot}, keeping the sign conventions of
\eqref{eq:logkernel} and \eqref{eq:stokeskernels}. The pressure of the Stokes single layer is
\begin{equation}\label{eq:pressureSdirect}
 p^{S}(y)=-\frac{1}{2\pi}\int_\Gamma
 \frac{\boldsymbol r\cdot f(x)}{r^2}\dd s_x
 =-\int_\Gamma \nabla_xG(\boldsymbol r)\cdot f(x)\dd s_x,
 \qquad \boldsymbol r=x-y.
\end{equation}
Decompose $f=(f\cdot n)n+(f\cdot\tau)\tau$; the same split is used for
surfaces in \cite[Sect.~2]{BT2026}, which treats the tangential part by a
subtraction instead. The normal part of
\eqref{eq:pressureSdirect} is $-D[f\cdot n]$. For the tangential part,
$\partial_\tau G=\partial_sG$ and periodic integration by parts give
$-\int_\Gamma\partial_sG(f\cdot\tau)\dd s
=S[(f\cdot\tau)_s]$. Hence
\begin{equation}\label{eq:pressureS}
 p^S=S[(f\cdot\tau)_s]-D[f\cdot n].
\end{equation}
To compute $p^S$, we replace $S$ and $D$ by their regularized versions.

For $y\notin\Gamma$, the pressure of the Stokes double layer is
\begin{equation}\label{eq:pressureDdirect}
 p^D(y)=2\mu\int_\Gamma q(x)\cdot[\nabla_x^2G(\boldsymbol r)]n(x)\dd s_x
 =\frac{\mu}{\pi}\int_\Gamma\left[
 \frac{q\cdot n}{r^2}-2\frac{(q\cdot\boldsymbol r)(n\cdot\boldsymbol r)}{r^4}
 \right]\dd s_x ,
\end{equation}
where $\nabla_x^2G$ denotes the Hessian matrix of $G$.
This integral is hypersingular for $y\in\Gamma$. We therefore reduce it to $S$
and $D$, whose limits on $\Gamma$ follow from the continuity of $S$ and from
\eqref{eq:jump}. Writing $q$ and its derivatives in the
$(\tau,n)$ frame and using
$\tau_s=-\kappa n$, $n_s=\kappa\tau$, and $\Delta G=0$ away from the
source, integration by parts along $\Gamma$ gives
\begin{equation}\label{eq:pressureDreduce}
 \int_\Gamma q\cdot(\nabla_x^2G)n\dd s
 =\int_\Gamma(q_s\cdot n)\,\partial_\tau G\dd s
 -\int_\Gamma(q_s\cdot\tau)\,\partial_nG\dd s.
\end{equation}
Integrating the first term by parts once more yields
\begin{equation}\label{eq:pressureD}
 p^D=-2\mu S[(q_s\cdot n)_s]-2\mu D[q_s\cdot\tau],
 \qquad
 (q_s\cdot n)_s=q_{ss}\cdot n+\kappa q_s\cdot\tau .
\end{equation}

For the sum of the Stokeslet and stresslet integrals, define
\begin{equation}\label{eq:pressuredensities}
 \phi=(f\cdot\tau)_s-2\mu(q_s\cdot n)_s,
 \qquad
 \psi=f\cdot n+2\mu q_s\cdot\tau .
\end{equation}
By \eqref{eq:pressureS} and \eqref{eq:pressureD}, the pressure is
\begin{equation}\label{eq:pressurecombined}
 p=S\phi-D\psi.
\end{equation}
The pressure can therefore be computed with the regularized kernels and
coefficients of Section~\ref{sec:e1ansatz}.

The regularized pressure is
\begin{equation}\label{eq:pressurereg}
 p_\delta=S_\delta^{(p)}\phi-D_\delta^{(p)}\psi,
\end{equation}
where $D_\delta^{(p)}\psi$ is the subtracted double layer \eqref{eq:Dreg} with
density $\psi$.
For a periodic parametrization $x(t)$, we compute $q_s$, $(f\cdot\tau)_s$,
and $(q_s\cdot n)_s$ by Fourier differentiation, using
$\partial_s=|x'(t)|^{-1}\partial_t$.
Thus the hypersingular pressure kernel is replaced by tangential
derivatives of smooth periodic densities.

\subsection{Error estimates}\label{sec:stokesbounds}
The corrected Stokeslet and stresslet integrals have the same accuracy as
the potentials $S$ and $D$ in Proposition~\ref{prop:lap}; the proofs are given in
Section~\ref{sec:proofs}.
\begin{proposition}[Uniform Stokes accuracy]\label{prop:stokes}
Let $\Gamma$, $\Lambda$, and $\delta_0$ be as in
Proposition~\ref{prop:lap}, let $p=3$, $5$, or $7$, and let the vector
densities $f$ and $q$ be smooth. Let $\SL_{\delta,{\rm corr}}^{(p)}f$ and
$\DL_{\delta,{\rm corr}}^{(p)}q$ be given by \eqref{eq:SLcorr} and
\eqref{eq:DLcorr}, with the kernels \eqref{eq:Udelta} and \eqref{eq:Hdelta},
$\mathcal{N}^{(p)}$ and $\mathcal{V}$ from \eqref{eq:Lclosed} and \eqref{eq:C2}, and the
coefficients \eqref{eq:c3}--\eqref{eq:c7} at $\lambda=b/\delta$, and let
$\SL f$ and $\DL q$ be the potentials \eqref{eq:stokeslayers}.
Then for $0<\delta\le\delta_0$,
\begin{equation}\label{eq:stokesbound}
 \sup_{x_0\in\Gamma,\ |b|\le\Lambda\delta}
 \Bigl(|\SL_{\delta,{\rm corr}}^{(p)}f-\SL f|
 +|\DL_{\delta,{\rm corr}}^{(p)}q-\DL q|\Bigr)\le C\delta^{p}.
\end{equation}
The constant $C$
depends on finitely many derivatives of $\Gamma$ and the densities, and on
$p$ and $\Lambda$, but not on $\delta$, $b$, or $x_0$.
\end{proposition}
Since the pressure is computed from $S$ and $D$, the same holds for
the pressure.
\begin{proposition}[Uniform pressure accuracy]\label{prop:pressure}
Let $\Gamma$, $\Lambda$, $\delta_0$, and the order be as in
Proposition~\ref{prop:lap}, and let $f$ and $q$ be smooth, so that the
pressure densities $\phi$ and $\psi$ of \eqref{eq:pressuredensities} are
smooth. Let $p_\delta$ be the regularized pressure \eqref{eq:pressurereg},
with $S_\delta^{(p)}$ and $D_\delta^{(p)}$ given by \eqref{eq:Sreg} and
\eqref{eq:Dreg}. Then for $0<\delta\le\delta_0$,
\begin{equation}\label{eq:pressurebound}
 \sup_{x_0\in\Gamma,\ |b|\le\Lambda\delta}
 |p_\delta-p|\le C\delta^{p}.
\end{equation}
\end{proposition}

\begin{remark}\label{rem:oncurvestokes}
As for the Laplace double layer (Remark~\ref{rem:oncurve}), for $y\in\Gamma$
the stresslet kernel contracted with the normal, $T_{ijk}(x-y)n_k(x)$, is a
smooth function of $x\in\Gamma$, with limit
$\kappa(y)\tau_i(y)\tau_j(y)/(2\pi)$ as $x\to y$ \cite{GKM1996}. With this
diagonal value, the stresslet integral needs no regularization on the curve,
whereas the Stokeslet and the single-layer part of the pressure remain singular
on $\Gamma$.
\end{remark}

\section{Error analysis}\label{sec:analysis}
We now derive the kernel representations and local corrections of Sections~\ref{sec:e1ansatz}
and~\ref{sec:stokes}, establishing the uniform error bounds of Propositions~\ref{prop:lap},
\ref{prop:stokes}, and~\ref{prop:pressure}. The analysis constructs a local asymptotic expansion
in a tubular neighborhood of the closest point $x_0$. Following \cite{BL2001,BT2024}, each error term
factors into a geometric Taylor coefficient and an error moment depending solely on $\lambda=b/\delta$.

\subsection{The single layer}\label{sec:baserror}
Let $y=x_0+bn_0$ be a target as in \eqref{eq:offset}, so that $b$ is the
signed normal distance from the curve to $y$, and let $\lambda=b/\delta$. Under
the hypotheses of Proposition~\ref{prop:lap}, the closest point $x_0$ is
unique and $1+b\kappa_0\ge\frac12$, where $\kappa_0$ is the curvature at $x_0$.
Let $\sigma$ be signed
arclength from $x_0$, positive in the direction of $\tau_0$, and let
$\boldsymbol r(\sigma)=x(\sigma)-y=x(\sigma)-x_0-bn_0$. Taylor expansion of the
curve about $x_0$, with $\tau_s=-\kappa n$ and $n_s=\kappa\tau$, gives
\begin{equation}\label{eq:localgeom}
 \boldsymbol r=-bn_0+\sigma\tau_0-\tfrac12\kappa_0\sigma^2n_0
 -\tfrac16\sigma^3(\kappa_0^2\tau_0+\kappa_sn_0)+O(\sigma^4),
\end{equation}
where $\kappa_s$ is the arclength derivative of the curvature at $x_0$.
Hence $\boldsymbol r\cdot n(\sigma)=-b+O(\sigma^2)$ and
$r^2-b^2=\alpha^2\sigma^2+\tfrac13b\kappa_s\sigma^3+O(\sigma^4)$, with
$\alpha=\sqrt{1+b\kappa_0}$. Since $\alpha^2\ge\frac12$, the coordinate
$\eta=\operatorname{sgn}(\sigma)\sqrt{r^2-b^2}$ along the curve is a smooth,
invertible function of $\sigma$ near $0$, uniformly in $x_0$ and $b$, and
$r^2=\eta^2+b^2$. Taking the square root and inverting the series gives
\begin{equation}\label{eq:eta}
 \eta=\alpha\sigma+\frac{b\kappa_s}{6\alpha}\sigma^2+O(\sigma^3),\qquad
 \sigma=\frac{\eta}{\alpha}-\frac{b\kappa_s}{6\alpha^4}\eta^2+O(\eta^3),
\end{equation}
so the Jacobian $\mathcal{J}(\eta)=\dd\sigma/\dd\eta$ satisfies
$\mathcal{J}(0)=1/\alpha$. The coordinate $\eta$ is the two-dimensional
analogue of the one used in \cite[Sect.~3]{BT2024}. In this coordinate, a radial kernel becomes a function of
$\eta^2+b^2$, and the density, the Jacobian, and the numerators of the Stokes
kernels have Taylor expansions in $\eta$.
With $\xi=(\eta^2+b^2)/\delta^2$, the Gaussian $e^{-\xi}$ is concentrated where
$\eta^2+b^2=O(\delta^2)$, that is, on $|\eta|=O(\delta)$ when $b=O(\delta)$,
uniformly for $|\lambda|\le\Lambda$. The part of the curve outside a fixed
neighborhood of $x_0$ therefore contributes only exponentially small terms, and
only the behavior of the integrand near $\eta=0$ matters.

For the Laplace single layer, the kernel difference is
$G_\delta^{(1)}-G=\frac{1}{4\pi}E_1(\xi)$. We combine the density and the
Jacobian into the smooth function
\begin{equation}\label{eq:SLamp}
 w(\eta;b)=f(x(\sigma(\eta)))\,\mathcal{J}(\eta),
\end{equation}
and insert a smooth even cutoff function that equals one near $\eta=0$ and
vanishes outside the coordinate patch. Suppressing the cutoff in the
notation and extending the integral to $\R$, which changes it by an
exponentially small amount, we obtain
\begin{equation}\label{eq:baseerror}
 S_\delta^{(1)}f-Sf
 =\frac{1}{4\pi}\int_\R w(\eta;b)E_1
 \Bigl(\frac{\eta^2+b^2}{\delta^2}\Bigr)\dd\eta
 +\text{exponentially small terms}.
\end{equation}
Write the even part of $w$ as
\begin{equation}\label{eq:aeven}
 w_{\rm ev}(\eta;b)=\tfrac12\bigl[w(\eta;b)+w(-\eta;b)\bigr]
 =\sum_{\substack{n\ge0\\ n\ {\rm even}}}w_n(x_0,b)\eta^n;
\end{equation}
finite Taylor expansions with remainder suffice throughout, so analyticity is
not required. Only $w_{\rm ev}$ contributes, because the kernel factor is
even in $\eta$ for fixed $b$, so the odd part of the integrand integrates to
zero.
Substituting \eqref{eq:aeven} into \eqref{eq:baseerror} and setting
$\eta=\delta t$ gives \eqref{eq:baseexp}. For the kernel \eqref{eq:Greg},
the same steps give \eqref{eq:SLerror}.

\subsection{The coefficients}\label{sec:coeff}
We first derive the form of the correction terms. Differentiating
$E_1(\xi)=\int_1^\infty e^{-\xi v}v^{-1}\dd v$ $k$ times and substituting
$u=\xi v$ gives
\begin{equation}\label{eq:derivativeidentity}
 \xi^kE_1^{[k]}(\xi)=(-1)^k\int_\xi^\infty u^{k-1}e^{-u}\dd u
 =(-1)^k(k-1)!\,e^{-\xi}\sum_{\ell=0}^{k-1}\frac{\xi^\ell}{\ell!},
 \qquad k\ge1,
\end{equation}
where the second equality follows by integrating by parts $k-1$ times.
Differentiating the integral gives
$\frac{\dd}{\dd\xi}[\xi^kE_1^{[k]}(\xi)]=(-1)^{k+1}\xi^{k-1}e^{-\xi}$.

We now solve the moment conditions \eqref{eq:M}. The polynomial factor of
$\xi^kE_1^{[k]}$ has degree $k-1$, which increases with $k$, so the
conditions can be solved by successive elimination, as we now show.
The $E_1$ term in \eqref{eq:F} contributes a function
$\mathcal{H}_n(\lambda)$, $n=0,2,4,\ldots$, to each moment, and each
correction term $\xi^kE_1^{[k]}$ contributes an elementary Gaussian integral.
Together they give a small linear system for the $a_{k}^{(p)}$. For even
$n\ge0$, inserting the integral definition of $E_1$ and interchanging the
order of integration gives
\begin{equation}\label{eq:basemoment}
 \begin{aligned}
 \int_\R t^nE_1(t^2+\lambda^2)\dd t
 &=\int_\R t^n\!\int_1^\infty
   \frac{e^{-(t^2+\lambda^2)v}}{v}\dd v\,\dd t\\
 &=\Gamma\!\left(\frac{n+1}{2}\right)
   \int_1^\infty v^{-(n+3)/2}e^{-\lambda^2v}\dd v\\
 &=\Gamma\!\left(\frac{n+1}{2}\right)
   e^{-\lambda^2}\mathcal{H}_n(\lambda).
 \end{aligned}
\end{equation}
Here $\mathcal{H}_n$ is defined in \eqref{eq:Hdef}; the factor
$e^{-\lambda^2}$ is taken out because it is common to all terms of the
moment conditions.

By \eqref{eq:derivativeidentity}, the $k$-th correction term in
\eqref{eq:F} is a polynomial in $\xi=t^2+\lambda^2$ times
$e^{-(t^2+\lambda^2)}$. Expanding in powers of $t^2$ and applying
$\int_\R t^{2i}e^{-t^2}\dd t=\Gamma(i+\tfrac12)$ to each term gives, for even
$n\ge0$,
\begin{equation}\label{eq:polyterm}
 \int_\R t^n\xi^kE_1^{[k]}(\xi)\dd t
 =(-1)^k(k-1)!\,e^{-\lambda^2}\sum_{\ell=0}^{k-1}\frac{1}{\ell!}
 \sum_{i=0}^{\ell}\binom{\ell}{i}\lambda^{2(\ell-i)}
 \Gamma\!\left(i+\frac{n+1}{2}\right).
\end{equation}
After division by the common factor
$\Gamma((n+1)/2)e^{-\lambda^2}$, the $k$-th correction term contributes a
polynomial in $\lambda^2$ of degree $k-1$ with leading term
$(-1)^k\lambda^{2(k-1)}$, and the moment conditions become
\begin{equation}\label{eq:linsys}
 \mathcal{H}_n-a_{1}^{(p)}
 +\left(\lambda^2+\frac{n+3}{2}\right)a_{2}^{(p)}
 -\left[\lambda^4+(n+3)\lambda^2
 +\frac{(n+3)(n+5)}{4}\right]a_{3}^{(p)}=0,
\end{equation}
for $n=0,2,\ldots,p-3$. Here $a_{k}^{(p)}=0$ for $k>m=(p-1)/2$; the alternating
signs come from the factors $(-1)^k$ in $\xi^kE_1^{[k]}$.

Since the coefficient of $a_{k}^{(p)}$ in \eqref{eq:linsys} is a polynomial of
degree $k-1$ in $n$, taking differences of the equations for consecutive
values of $n$ eliminates the unknowns one at a time. For $p=7$, subtracting
the equation for $n=0$ from that for $n=2$ gives
$a_2^{(7)}=\mathcal{H}_0-\mathcal{H}_2+(2\lambda^2+5)a_3^{(7)}$, the second
difference gives $a_3^{(7)}=\tfrac12(\mathcal{H}_0-2\mathcal{H}_2+\mathcal{H}_4)$,
and the equation for $n=0$ then gives $a_1^{(7)}$. The cases $p=5$ and $p=3$
follow in the same way, which gives \eqref{eq:c3}--\eqref{eq:c7}.

\subsection{The double layer}\label{sec:dlanalysis}
For the local error analysis of the subtracted double layer
\eqref{eq:Dreg}, define the numerator of its integrand, without the factor
$d^{(p)}/(2\pi r^2)$ of \eqref{eq:Kreg} and with the Jacobian included,
\begin{equation}\label{eq:AD}
 w_D(\eta;b)=\bigl[n(x)\cdot(x-y)\bigr]\bigl[g(x)-g(x_0)\bigr]\mathcal{J}(\eta),
 \qquad w_D(0;b)=0,
\end{equation}
with the Jacobian $\mathcal{J}$ of Section~\ref{sec:baserror}. The terms
$\chi g(x_0)$ in \eqref{eq:Dreg} and in
$Dg=\int_\Gamma K\,(g-g(x_0))\dd s+\chi g(x_0)$ cancel, and
$K_\delta^{(p)}-K=-n(x)\cdot(x-y)\,Q^{(p)}(\xi;\lambda)/(2\pi r^2)$ by
\eqref{eq:Kreg}, since $d^{(p)}=1-Q^{(p)}$, with $Q^{(p)}$ as in
\eqref{eq:PQT}. The argument of
Section~\ref{sec:baserror}, applied to the even part of $w_D$, which begins
with $\eta^2$, then gives
\begin{equation}\label{eq:DLerror}
 D_\delta^{(p)}g-Dg\sim-\frac{1}{2\pi}
 \sum_{\substack{n\ge0\\ n\ {\rm even}}}w_{D,n+2}(x_0,b)\,
 \delta^{n+1}J_n^{(p)}(\lambda),
\end{equation}
with the moments $J_n^{(p)}$ of \eqref{eq:Jdef}, where $w_{D,n+2}(x_0,b)$ is the
coefficient of $\eta^{n+2}$ in the Taylor expansion of $w_D(\eta;b)$. The moments
$J_n^{(p)}$ are multiples of the single-layer moments.

\begin{lemma}[Relation between the moments]\label{lem:transfer}
For every even $n\ge0$ and every $\lambda$,
\begin{equation}\label{eq:transfer}
 J_n^{(p)}(\lambda)=\frac{n+1}{2}\,I_n^{(p)}(\lambda).
\end{equation}
\end{lemma}
\begin{proof}
Put $u=t^2+\lambda^2$. By \eqref{eq:PQT},
$Q^{(p)}(u;\lambda)=-u\partial_u s^{(p)}(u;\lambda)$. Integrating by parts and using
$\frac{\dd}{\dd t}s^{(p)}(u;\lambda)=2t\partial_u s^{(p)}(u;\lambda)$,
\[
 J_n^{(p)}=-\int_\R t^{n+2}\partial_u s^{(p)}(u;\lambda)\dd t
 =-\frac12\int_\R t^{n+1}\frac{\dd}{\dd t}s^{(p)}(u;\lambda)\dd t
 =\frac{n+1}{2}\int_\R t^ns^{(p)}(u;\lambda)\dd t .
\]
The boundary terms vanish at infinity because $s^{(p)}$ decays like a Gaussian
times algebraic factors. At $\lambda=0$, split the integration at $t=0$; the
additional endpoint contributions vanish since $t^{n+1}\log(t^2)\to0$.
\end{proof}

Hence the coefficients that cancel the first $m$ moments
$I_n^{(p)}$ also cancel the first $m$ moments $J_n^{(p)}$, and the double
layer has order $p$ with the same coefficients. This relation is the
two-dimensional analogue of the relation between single- and double-layer
moments in \cite[eq.~(34)]{BT2026}.

Singularity subtraction is needed for uniform accuracy near the curve.
Without it, the regularized double layer $\int K_\delta^{(p)}g\dd s$ does not
converge uniformly for $b=\lambda\delta$ with fixed $\lambda\ne0$, because
the regularized kernel does not satisfy Gauss's identity \eqref{eq:chi}. By
\eqref{eq:chi}--\eqref{eq:Kreg}, the kernel difference integrates to
\begin{equation}\label{eq:bias}
 \int_\Gamma\bigl(K_\delta^{(p)}-K\bigr)\dd s_x
 =\beta^{(p)}(\lambda)\bigl(1+O(\delta)\bigr),\qquad
 \beta^{(p)}(\lambda)=\frac{\lambda}{2\pi}\int_\R
 \frac{Q^{(p)}(t^2+\lambda^2;\lambda)}{t^2+\lambda^2}\dd t ,
\end{equation}
so without subtraction the error tends to
$g(x_0)\beta^{(p)}(\lambda)$ as $\delta\to0$, which is
$O(1)$. Singularity subtraction removes this error exactly. Although
$Dg=\int K(g-\bar g)\dd s+\chi\bar g$ holds for any constant $\bar g$, we subtract
the value at the closest point, since then the numerator vanishes at
$\eta=0$ and the error expansion \eqref{eq:DLerror} holds.

\subsection{The Stokeslet}\label{sec:stokeslet}

Substituting \eqref{eq:Greg} into \eqref{eq:Udelta} and comparing with
\eqref{eq:stokeskernels} gives
\begin{equation}\label{eq:Udefect}
 S_{ij,\delta}^{(p)}-S_{ij}
 =-\frac{\delta_{ij}}{8\pi\mu}s^{(p)}(\xi;\lambda)
 -\frac{1}{4\pi\mu}\frac{r_ir_j}{r^2}Q^{(p)}(\xi;\lambda).
\end{equation}
The first term is a multiple of the regularizing function, whose
moments vanish by \eqref{eq:M}. The second term contains $r_ir_j/r^2$. By
\eqref{eq:localgeom} and \eqref{eq:eta}, the numerator $r_ir_j$ has the expansion
\begin{equation}\label{eq:tensorexp}
 r_ir_j=b^2(n_0)_i(n_0)_j
 -\frac{b}{\alpha}\eta\bigl[(n_0)_i(\tau_0)_j+(\tau_0)_i(n_0)_j\bigr]
 +\eta^2\Theta_{ij}(b)+O(\eta^3),
\end{equation}
where $\Theta_{ij}(b)$ is a smooth tensor coefficient, bounded for $b$
in the band around the curve, with
$\Theta_{ij}(0)=(\tau_0)_i(\tau_0)_j$. After division by
$r^2=\eta^2+b^2$ and the scaling $\eta=\delta t$, an even term $\eta^{n+2}$,
$n\ge0$, gives $\delta^{n+1}J_n^{(p)}$, with $J_n^{(p)}$ as in
\eqref{eq:Jdef}. These are multiples of the $I_n^{(p)}$ by
Lemma~\ref{lem:transfer}, so they impose no new conditions. The only even
term without a factor $\eta^2$ is $b^2(n_0)_i(n_0)_j$, and it gives
$\delta^{-1}b^2\int_\R Q^{(p)}(t^2+\lambda^2;\lambda)(t^2+\lambda^2)^{-1}\dd t
=|b|\,\mathcal{N}^{(p)}(\lambda)$, where
\begin{equation}\label{eq:Ndef}
 \mathcal{N}^{(p)}(\lambda)=|\lambda|\int_\R\frac{Q^{(p)}(t^2+\lambda^2;\lambda)}{t^2+\lambda^2}\dd t,
 \qquad \lambda\ne0.
\end{equation}
The conditions \eqref{eq:M} do not cancel this term, and we remove it by a
local correction.

The integral in \eqref{eq:Ndef} behaves like $\pi/|\lambda|$ as $\lambda\to0$,
but $\mathcal{N}^{(p)}$ has the closed form \eqref{eq:Lclosed}, which is
bounded and equal to $\pi$ at $\lambda=0$. To see this, write
$P^{(p)}(\xi;\lambda)=1+\xi Z$, where
$Z=\sum_{k=1}^m(-1)^ka_k^{(p)}(\lambda)\xi^{k-1}$, and split the integrand
into two parts. The first part, $e^{-\xi}/(t^2+\lambda^2)$, gives
$\pi\operatorname{erfc}|\lambda|$, because $\int_\R e^{-t^2}/(t^2+\lambda^2)\dd t
=(\pi/|\lambda|)e^{\lambda^2}\operatorname{erfc}|\lambda|$. The second part,
$\xi Z(\xi)e^{-\xi}/(t^2+\lambda^2)=Z(\xi)e^{-\xi}$, is a polynomial times a
Gaussian; integrating it term by term as in \eqref{eq:polyterm} gives the
remaining terms of \eqref{eq:Lclosed}.

To remove the remaining Stokeslet term, we need the constant term of the
even part of the numerator, including the Jacobian $\mathcal{J}$, and set
\begin{equation}\label{eq:ADU}
 w^U_i(\eta;b)=r_ir_j(\eta;b)\,f_j(x(\sigma(\eta)))\,\mathcal{J}(\eta).
\end{equation}
The second term of \eqref{eq:Udefect} then contributes
$-\frac{1}{4\pi\mu}\int_\R w^U_i(\eta;b)Q^{(p)}(\xi;\lambda)/(\eta^2+b^2)\dd\eta$
to the error. Since $r_ir_jf_j=b^2(n_0)_i(n_0\cdot f_0)+O(\eta)$ by
\eqref{eq:tensorexp} and $\mathcal{J}(0)=\alpha^{-1}$ by \eqref{eq:eta}, the
constant term of $w^U_i$ is $(b^2/\alpha)(n_0)_i(n_0\cdot f_0)$, and by
\eqref{eq:Ndef} its contribution to the regularization error is
$-\frac{|b|\mathcal{N}^{(p)}(\lambda)}{4\pi\mu\alpha}(n_0)_i(n_0\cdot f_0)$, which
vanishes as $b\to0$. Subtracting it gives the correction in
\eqref{eq:SLcorr}.

\subsection{The stresslet}\label{sec:stresslet}

From \eqref{eq:dfactors},
$\tilde d^{(p)}=1-\xi^2\partial_\xi^2s^{(p)}$ at fixed $\lambda$.
Comparing \eqref{eq:Hdelta} with \eqref{eq:stokeskernels} therefore gives the
exact identity
\begin{equation}\label{eq:Hdefect}
 T_{ijk,\delta}^{(p)}-T_{ijk}
 =-\frac{1}{\pi}\frac{r_ir_jr_k}{r^4}
 R^{(p)}(\xi;\lambda),
 \qquad R^{(p)}(\xi;\lambda)=1-\tilde d^{(p)}(\xi;\lambda)
 =\xi^2\partial_\xi^2s^{(p)}(\xi;\lambda) .
\end{equation}
We subtract the density value at $x_0$ as in \eqref{eq:Dreg} and
include the Jacobian in the numerator, which after contraction with
$n_k(x)$ is
\begin{equation}\label{eq:ADH}
 w^H_i(\eta;b)=r_ir_j(\boldsymbol r\cdot n(x))\,(q_j(x)-q_{0,j})\,
 \mathcal{J}(\eta).
\end{equation}
The kernel difference \eqref{eq:Hdefect} then contributes
$-\frac{1}{\pi}\int w^H_i(\eta;b)R^{(p)}(\xi;\lambda)/(\eta^2+b^2)^2\dd\eta$,
and the stresslet moments are
\begin{equation}\label{eq:Lfamily}
 K_n^{(p)}(\lambda)=\int_\R\frac{t^{n+2}}{(t^2+\lambda^2)^2}
 R^{(p)}(t^2+\lambda^2;\lambda)\dd t,
 \qquad n=0,2,4,\ldots.
\end{equation}
Integrating by parts twice reduces them to the scalar moments. With
$u=t^2+\lambda^2$ and
$\partial_u^2s^{(p)}(u;\lambda)=\frac{1}{2t}\frac{\dd}{\dd t}
\partial_us^{(p)}(u;\lambda)$,
\[
 K_n^{(p)}=\int_\R t^{n+2}\partial_u^2s^{(p)}(u;\lambda)\dd t
 =-\frac{n+1}{2}\int_\R t^n\partial_us^{(p)}(u;\lambda)\dd t
 =\frac{(n+1)(n-1)}{4}\int_\R t^{n-2}s^{(p)}(u;\lambda)\dd t ,
\]
so that
\begin{equation}\label{eq:Ltransfer}
 K_n^{(p)}=\frac{(n+1)(n-1)}{4}I_{n-2}^{(p)}
 \quad(n=2,4,6,\ldots),\qquad
 |\lambda|K_0^{(p)}=\tfrac12\mathcal{N}^{(p)}\quad(\lambda\ne0).
\end{equation}
In particular,
$K_2^{(p)}=\tfrac34I_0^{(p)}$ and
$K_4^{(p)}=\tfrac{15}{4}I_2^{(p)}$. The relation for $K_0^{(p)}$ follows from a
single integration by parts:
$K_0^{(p)}=\int_\R t^2\partial_u^2s^{(p)}(u;\lambda)\dd t
=-\tfrac12\int_\R\partial_u s^{(p)}(u;\lambda)\dd t
=\tfrac12\int_\R Q^{(p)}(u;\lambda)u^{-1}\dd t$ by \eqref{eq:PQT}, which
is $\mathcal{N}^{(p)}/(2|\lambda|)$ by \eqref{eq:Ndef}. Thus the
conditions \eqref{eq:M} cancel all stresslet moments except $K_0^{(p)}$.

Since $q-q_0$ vanishes at $\eta=0$, the even part of $w^H$ begins with a
multiple of $\eta^2$, and this term pairs with $K_0^{(p)}$. To compute its
coefficient, we expand the numerator through second order in $\sigma$, before multiplying by
$\mathcal{J}$:
\begin{equation}\label{eq:Cexp}
 -b^3n_0(n_0\cdot q_s)\sigma
 +\Bigl\{b^2\bigl[n_0(\tau_0\cdot q_s)+\tau_0(n_0\cdot q_s)\bigr]
 -\tfrac12b^3n_0(n_0\cdot q_{ss})\Bigr\}\sigma^2+O(\sigma^3),
\end{equation}
using \eqref{eq:localgeom}. The first term is odd in $\sigma$ and the second
is even, but after the change of variable the first term also contributes
to the even part. Substituting
the series for $\sigma$ in \eqref{eq:eta} and its
derivative into \eqref{eq:Cexp} and collecting the coefficient of $\eta^2$
gives $b^2\mathcal{V}$, with the vector $\mathcal{V}$ of \eqref{eq:C2}.
The last term of $\mathcal{V}$ comes from the $\eta^2$ term of $\sigma$ acting
on the first term of \eqref{eq:Cexp}. After $\eta=\delta t$, the term
$b^2\mathcal{V}\eta^2$ contributes
$-\pi^{-1}\delta^{-1}b^2K_0^{(p)}\mathcal{V}=-|b|\mathcal{N}^{(p)}(\lambda)\mathcal{V}/(2\pi)$
to the error, by \eqref{eq:Ltransfer}. Subtracting it gives the
correction in \eqref{eq:DLcorr}, which thus involves the same function
$\mathcal{N}^{(p)}$ as the Stokeslet correction \eqref{eq:SLcorr}.

\subsection{Proofs of the error estimates}\label{sec:proofs}
We now prove Propositions~\ref{prop:lap}, \ref{prop:stokes},
and~\ref{prop:pressure}.
\begin{proof}[Proof of Proposition~\ref{prop:lap}]
Expand $w_{\rm ev}$ through degree $p-3$ and $w_{D,{\rm ev}}$ through degree
$p-1$, with remainder. After the integrals are extended to $\R$, which
changes them by an exponentially small amount, the terms of the Taylor
polynomials integrate to zero by \eqref{eq:M} and Lemma~\ref{lem:transfer};
the coefficients \eqref{eq:c3}--\eqref{eq:c7} satisfy \eqref{eq:M} by
construction (Section~\ref{sec:coeff}).
The remainders are bounded uniformly by $C|\eta|^{p-1}$ for the single
layer and by $C|\eta|^{p+1}$ for the double layer, whose kernel contributes
the factor $(\eta^2+b^2)^{-1}$. After $\eta=\delta t$, both carry the
prefactor $\delta^{p}$, multiplying
\[
 \int_\R |t|^{p-1}\,|s^{(p)}(t^2+\lambda^2;\lambda)|\,\dd t
 \quad\text{and}\quad
 \int_\R \frac{|t|^{p+1}}{t^2+\lambda^2}\,|Q^{(p)}(t^2+\lambda^2;\lambda)|\,\dd t
 \le\int_\R |t|^{p-1}\,|Q^{(p)}(t^2+\lambda^2;\lambda)|\,\dd t,
\]
and these integrals are bounded uniformly for $|\lambda|\le\Lambda$: $s^{(p)}$ has
at most a logarithmic singularity, and $Q^{(p)}$ is a Gaussian times a polynomial
with coefficients bounded on $[-\Lambda,\Lambda]$, as noted in
Section~\ref{sec:props}. The part of the curve outside the
coordinate patch contributes an exponentially small amount.
\end{proof}

\begin{proof}[Proof of Proposition~\ref{prop:stokes}]
We localize the integrals as in the proof of Proposition~\ref{prop:lap}. The
terms involving $s^{(p)}$ vanish by \eqref{eq:M}. The other Stokeslet
moments are the $J_n^{(p)}$, and the stresslet moments $K_n^{(p)}$ with
$n\ge2$ are multiples of the $I_{n-2}^{(p)}$ by \eqref{eq:Ltransfer}; these
also vanish. The terms with $\mathcal{N}^{(p)}$ and $K_0^{(p)}$ are removed by
the local corrections in \eqref{eq:SLcorr} and \eqref{eq:DLcorr}. The Taylor
remainders are bounded by $C|\eta|^{p-1}$, $C|\eta|^{p+1}$, and
$C|\eta|^{p+3}$ in the scalar, Stokeslet tensor, and stresslet tensor
families, respectively, and $\eta=\delta t$ leaves
$\delta^{p}$ times integrals bounded uniformly in $\lambda$, since
$|t|^{p+1}/(t^2+\lambda^2)\le|t|^{p-1}$ and
$|t|^{p+3}/(t^2+\lambda^2)^2\le|t|^{p-1}$. Finally, $\mathcal{N}^{(p)}$ and $\mathcal{V}$ are bounded, so the corrections are
$O(|b|)$, vanish at $b=0$, and the estimate also holds on the curve.
\end{proof}

\begin{proof}[Proof of Proposition~\ref{prop:pressure}]
By \eqref{eq:pressurecombined} and \eqref{eq:pressurereg},
$p_\delta-p=(S_\delta^{(p)}-S)\phi
-(D_\delta^{(p)}-D)\psi$. Proposition~\ref{prop:lap} applies to the two
terms directly.
\end{proof}

\subsection{The discretization error}\label{sec:discrete-error}
The regularized integrals are computed with the periodic trapezoidal rule,
whose error is small when $\delta/h$ is large. Following
\cite{BT2024,BT2026}, we bound this error to see how $\delta$ should depend
on $h$.
Let $x(t)$ be a smooth $2\pi$-periodic parametrization of $\Gamma$ with
$\min_t|x'(t)|>0$, and let the densities be smooth. Let the targets lie in a
band around $\Gamma$ in which closest points are unique, with
$|b|\le\Lambda\delta$ and $1+b\kappa_0$ bounded below by a positive
constant. Let $F_{\delta,y}(t)$ be any scalar component of an integrand, and
$h=2\pi/N_\Gamma$. For each integer $\ell\ge2$,
\[
 \Bigl|h\sum_{j=0}^{N_\Gamma-1}F_{\delta,y}(jh)
       -\int_0^{2\pi}F_{\delta,y}(t)\dd t\Bigr|
 \le C_\ell h^\ell\|\partial_t^\ell F_{\delta,y}\|_{L^1(0,2\pi)},
\]
since the error of the rule is $2\pi$ times the sum of the Fourier
coefficients of $F_{\delta,y}$ at the nonzero multiples of $N_\Gamma$, and
integration by parts bounds the coefficient at $k$ by
$|k|^{-\ell}\|\partial_t^\ell F_{\delta,y}\|_{L^1}/(2\pi)$.

A crude bound on the derivatives is
\[
 \|\partial_t^\ell F_{\delta,y}\|_{L^1}
 \le C_\ell(1+|\log\delta|)\delta^{-\ell-1},\qquad 0<\delta\le1,
\]
uniformly in $x_0$ and $|\lambda|\le\Lambda$. With the coefficients held
fixed in $t$, the single-layer and Stokeslet kernels are smooth functions of
$(x-y)/\delta$ plus a multiple of $\log\delta$, so they are
$O(1+|\log\delta|)$ on $\Gamma$, and each $t$-derivative brings a factor
$\delta^{-1}$. The double-layer and stresslet kernels are $\delta^{-1}$ times
smooth bounded functions of $(x-y)/\delta$, since $d^{(p)}=O(\xi)$ and
$\tilde d^{(p)}=O(\xi^2)$ at $\xi=0$ and $1-d^{(p)}$ and $1-\tilde d^{(p)}$
decay exponentially. The chain and product rules give the bound. The
pressure uses the same kernels, with the densities $\phi$ and $\psi$.

Let $V$ be any of the potentials and $V_{\delta,h}$ its trapezoidal
approximation, including the local corrections and the term $\chi g_0$,
which involve no quadrature. With exact data and arithmetic,
Propositions~\ref{prop:lap}, \ref{prop:stokes}, and~\ref{prop:pressure} and
the two bounds above give
\[
 |V_{\delta,h}-V|
 \le C\delta^{p}+C_\ell h^\ell\delta^{-\ell-1}(1+|\log\delta|).
\]
For $\delta=ch^q$ with $0<q<1$, the second term is
$O(h^{\ell(1-q)-q}|\log h|)$, which is $o(h^{pq})$ if
\begin{equation}\label{eq:quadrature-sufficient-r}
 \ell(1-q)>q(p+1).
\end{equation}
Since $\ell$ is arbitrary for smooth data, the error is $O(h^{pq})$ for
every $q<1$. When $\delta$ is proportional to $h$, the second term is
$O(h^{-1}|\log h|)$ for every $\ell$, and the bound does not imply
convergence.

\section{Numerical examples}\label{sec:num}
We test the regularized kernels on the unit circle, a kite, and a
starfish (the kite and starfish are shown in Figure~\ref{fig:allgrid-maps}).
In this section, target points $y=x(t_0)+b\,n(t_0)$ are specified directly by the
closest-point parameter $t_0$ and signed normal distance $b$. On the Cartesian grids of
Section~\ref{sec:num-allgrid}, $t_0$ is found by Newton's method applied to
$(x(t)-y)\cdot x'(t)=0$, starting from the closest of $2^{12}$ sample points on $\Gamma$,
with $b=(y-x_0)\cdot n_0$ (cf.\ \cite[Sect.~4]{BL2001}). The curvature $\kappa$ and its
derivative $\kappa_s$ are evaluated from the parametrization of $\Gamma$.

Unless stated otherwise, the reported error is the maximum over the 97
parameters $t_0=2\pi(j+0.37)/97$ ($j=0,\dots,96$) and $b/\delta=-1/2,0,1/2$.
Targets on the curve ($b=0$) use the regularized kernels with $\lambda=0$
so that a single formula covers all targets, although on $\Gamma$ the
double-layer and stresslet kernels are smooth and need no regularization
(Remarks~\ref{rem:oncurve} and~\ref{rem:oncurvestokes}). Table~\ref{tab:rates} summarizes the convergence rates for the main tests.
All computations were done in double precision in Python with NumPy and SciPy \cite{SciPy2020}
on a MacBook Air with an Apple M3 chip and 8 GB of memory.

\subsection{Two choices of \texorpdfstring{$\delta$}{delta}}\label{sec:delta-choice}
We use two choices of $\delta$. With the parameter spacing $h=2\pi/N_\Gamma$,
the first is
\begin{equation}\label{eq:numscaling}
 \delta=c h^q,\qquad c=1,\qquad
 (p,q)=(3,2/3),(5,4/5),(7,5/7),
\end{equation}
so the regularization error is $O(h^{pq})$ and the predicted orders are
$pq=2,4,5$. Since $q<1$, $\delta/h\to\infty$, and
\eqref{eq:quadrature-sufficient-r} holds with $\ell=9,25,21$; for smooth curves and
densities the discretization error is therefore asymptotically smaller than the
regularization error. Throughout, we estimate the discretization error by the
change when $N_\Gamma$ is doubled with $\delta$ and the targets fixed.

The second choice, as in \cite{BT2026}, takes $\delta$ proportional to the mesh spacing:
$\delta=ch$ for the circle, and $\delta=c\Delta s$ for the kite and the starfish, whose
nodes are not equally spaced along the curve. Here
\begin{equation}\label{eq:numphysical}
 \Delta s
 :=\max_j|x(t_{j+1})-x(t_j)|,
 \qquad t_j=jh,
\end{equation}
is the largest distance between adjacent nodes, and rates are fitted against $\Delta s$.
The regularization error is $O(h^p)$, but the estimate of Section~\ref{sec:discrete-error}
does not imply convergence, so the discretization error need not become negligible. We
therefore fit rates only over ranges where the estimated discretization error is well below
the total error. For this choice we also maximize the error over closest points at nodes and
halfway between nodes, and for the circle also a quarter of the way.

\subsection{The Laplace equation}\label{sec:num-laplace}

\subsubsection{Single- and double-layer potentials \texorpdfstring{for a circle}{for a circle}}\label{sec:num-circle}
We first test the single and double layers separately, with $\Gamma$ the unit circle and
the density
\begin{equation}\label{eq:circle-density}
 f(t)=1+0.4\cos3t+0.3\sin5t,
 \qquad x(t)=(\cos t,\sin t).
\end{equation}
For $y=R(\cos\theta,\sin\theta)$ and $G=+\log r/(2\pi)$,
\begin{equation}\label{eq:circle-single-exact}
 Sf(y)=\log\max(R,1)-\frac{0.4}{6}\min(R^3,R^{-3})\cos3\theta
 -\frac{0.3}{10}\min(R^5,R^{-5})\sin5\theta,
\end{equation}
which is continuous at $R=1$. The double layer is
\begin{equation}\label{eq:circle-double-exact}
 Df(y)=\begin{cases}
 1+0.2R^3\cos3\theta+0.15R^5\sin5\theta,&R<1,\\
 -0.2R^{-3}\cos3\theta-0.15R^{-5}\sin5\theta,&R>1,
 \end{cases}
 \qquad \mathcal{K}f=\tfrac12.
\end{equation}
These formulas follow mode by mode from the Fourier expansion of
$\log|x-y|$ and satisfy $D^-f-D^+f=f$.

\begin{figure}[!htbp]\centering
 \begin{minipage}[t]{.49\linewidth}\centering
  \includegraphics[width=\linewidth]{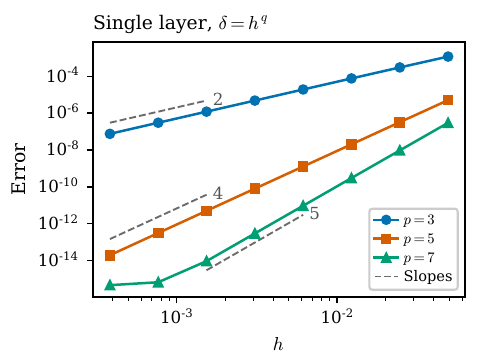}\par\vspace{-3pt}
  {\small (a)}
  \end{minipage}\hfill\begin{minipage}[t]{.49\linewidth}\centering
  \includegraphics[width=\linewidth]{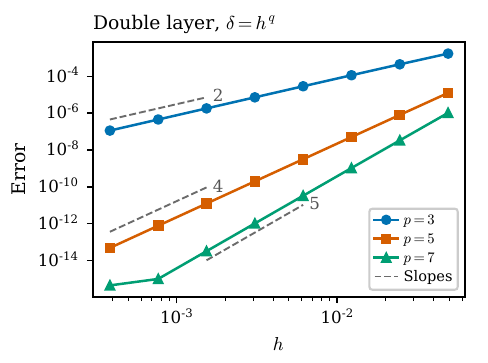}\par\vspace{-3pt}
  {\small (b)}
  \end{minipage}\par\smallskip
 \begin{minipage}[t]{.49\linewidth}\centering
  \includegraphics[width=\linewidth]{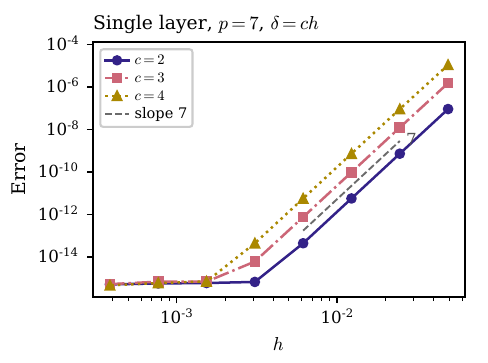}\par\vspace{-3pt}
  {\small (c)}
  \end{minipage}\hfill\begin{minipage}[t]{.49\linewidth}\centering
  \includegraphics[width=\linewidth]{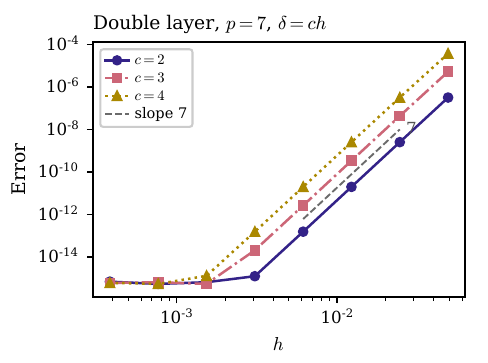}\par\vspace{-3pt}
  {\small (d)}
  \end{minipage}
 \caption{Single layer (a),(c) and double layer (b),(d) for the unit circle:
 (a),(b) $\delta=h^q$; (c),(d) $p=7$, $\delta=ch$.}
 \label{fig:circle-layers}
\end{figure}

For $p=3,5,7$, both
layers converge at the predicted orders $2,4,5$ (Table~\ref{tab:rates},
Figure~\ref{fig:circle-layers}(a),(b)).
With $\delta=c h$, $c=2,3,4$, both seventh-order layers follow the
reference slope seven until they reach roundoff
(Figure~\ref{fig:circle-layers}(c),(d)).

\begin{table}[tbp]\centering\small
\caption{Fitted convergence rates in $h$ (top) and in $\Delta s$ (bottom). The fits use
$N_\Gamma=2^{11}$--$2^{14}$ for the circle with $p=3,5$ and $2^{9}$--$2^{12}$ with
$p=7$; $2^{12}$--$2^{14}$ for the kite and the starfish with $p=3,5$ and
$2^{13}$--$2^{15}$ with $p=7$; and $2^{12}$--$2^{14}$ in the bottom block, except
$2^{11}$--$2^{13}$ for the starfish pressure.}
\label{tab:rates}
\begin{tabular}{llccc}
\toprule
Test & Quantity & $p=3$ & $p=5$ & $p=7$\\
\midrule
Predicted, $\delta=h^q$ & & 2 & 4 & 5\\
Circle & single layer & 2.00 & 4.00 & 5.00\\
Circle & double layer & 2.00 & 4.00 & 5.00\\
Kite & $Sf+Dg$ & 1.91 & 3.95 & 4.90\\
Starfish & velocity & 1.93 & 3.92 & 4.88\\
Starfish & pressure & 1.93 & 3.94 & 4.91\\
\midrule
 & & $c=2$ & $c=3$ & $c=4$\\
\midrule
Kite, $p=7$, $\delta=c \Delta s$ & $Sf+Dg$ & 6.94 & 6.88 & 6.80\\
Starfish, $p=7$, $\delta=c \Delta s$ & velocity & 6.94 & 6.87 & 6.77\\
Starfish, $p=7$, $\delta=c \Delta s$ & pressure & 6.80 & 6.62 & 6.38\\
\bottomrule
\end{tabular}
\end{table}

\subsubsection{A Laplace solution \texorpdfstring{for a kite}{for a kite}}
\label{sec:num-kite}
Next we test the sum $Sf+Dg$ of single and double layers, with densities on the kite-shaped
curve $\Gamma_L$,
\begin{equation}\label{eq:laplace-kite}
 x_L(t)=s_K(\cos t+0.65\cos2t-0.65,\ 1.5\sin t),
 \qquad s_K=\frac{179}{300}.
\end{equation}
With this scaling, $\max_t\|x_L(t)\|_\infty=0.895$. Let $z=x_1+\mathrm i x_2$ and set
$v^\mp=\operatorname{Re}F^\mp(z)$, where
\begin{align}
 F^-(z)={}&0.12+0.20z+(0.07-0.04\mathrm i)z^2
 +\frac{0.10}{z-\zeta_1}-\frac{0.08\mathrm i}{z-\zeta_2},\label{eq:laplace-Fminus}\\
 F^+(z)={}&\frac{0.16+0.03\mathrm i}{z-z_1}
 -\frac{0.07-0.04\mathrm i}{z-z_2},
 \label{eq:laplace-Fplus}
\end{align}
with
\begin{align*}
 \zeta_1&=2.40+1.10\mathrm i, & \zeta_2&=-2.50+0.80\mathrm i,\\
 z_1&=0.08+0.10\mathrm i, & z_2&=-0.22-0.08\mathrm i.
\end{align*}
The poles of $F^-$ lie outside the kite and those of $F^+$ lie inside, so
$v^-$ is harmonic inside and $v^+$ is harmonic outside. Since
$\nabla\operatorname{Re}F=(\operatorname{Re}F',-\operatorname{Im}F')$, the
normal derivatives are known exactly. With
\begin{equation}\label{eq:laplace-kite-densities}
 f=\partial_n v^+-\partial_n v^-,\qquad g=v^--v^+,
\end{equation}
the jump relations give $Sf+Dg=v^-$ inside, $v^+$ outside, and the average
$(v^-+v^+)/2$ on $\Gamma_L$.

\begin{figure}[!htbp]\centering
 \begin{minipage}[t]{.49\linewidth}\centering
  \includegraphics[width=\linewidth]{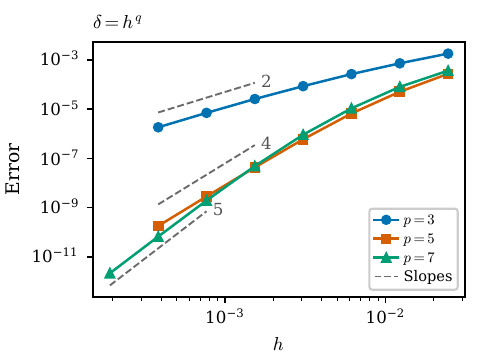}\par\vspace{-3pt}
  {\small (a)}
  \end{minipage}\hfill\begin{minipage}[t]{.49\linewidth}\centering
  \includegraphics[width=\linewidth]{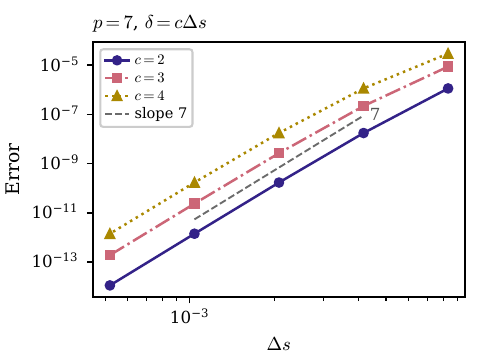}\par\vspace{-3pt}
  {\small (b)}
  \end{minipage}\par\smallskip
 \begin{minipage}[t]{.49\linewidth}\centering
  \includegraphics[width=\linewidth]{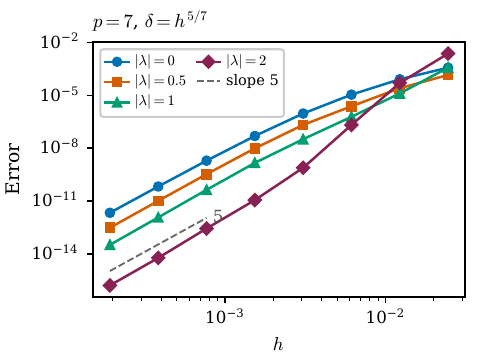}\par\vspace{-3pt}
  {\small (c)}
  \end{minipage}
 \caption{Laplace solution for the kite: (a) $\delta=h^q$; (b) $p=7$,
 $\delta=c\Delta s$; (c) $p=7$, $\delta=h^{5/7}$, for several $|\lambda|$.}
 \label{fig:laplace-kite-results}
\end{figure}

With $\delta$ as in \eqref{eq:numscaling}, the fitted rates are close
to the predicted orders, and with the seventh-order kernel and
$\delta=c \Delta s$, $c=2,3,4$, they are close to seven (Table~\ref{tab:rates},
Figure~\ref{fig:laplace-kite-results}(a),(b)).
Figure~\ref{fig:laplace-kite-results}(c) shows the seventh-order errors
separately for $|\lambda|=0,0.5,1,2$; for $\lambda\ne0$ we report the larger of
the errors on the two sides. On the finer meshes, the errors decrease as
$|\lambda|$ increases; at $|\lambda|=2$ the error on the finest mesh is at roundoff.

\subsection{Stokes flow}\label{sec:num-stokes}

\subsubsection{A Stokes solution \texorpdfstring{for a starfish}{for a starfish}}\label{sec:num-starfish}
For the Stokes test, we compute the velocity from the sum of the corrected
single and double layers, and the pressure from \eqref{eq:pressurereg}, for
the starfish $\Gamma_S$,
\begin{equation}\label{eq:stokes-starfish}
 x_S(t)=s_S(1+0.3\cos5t)(\cos t,\sin t),
 \qquad s_S=\frac{179}{260}.
\end{equation}
Again $\max_t\|x_S(t)\|_\infty=0.895$.
For $\mu=1$, define the interior stream function
\begin{equation}\label{eq:starfish-streamfunction}
 \Psi^-=\operatorname{Re}\{\bar zF(z)+W(z)\},\qquad
 u^-=(\Psi^-_{x_2},-\Psi^-_{x_1}),\qquad p^-=-4\operatorname{Im}F'(z),
\end{equation}
with $u^+=0$ and $p^+=0$.  The analytic functions are
\begin{align}
 F(z)={}&0.04+(0.055+0.025\mathrm i)z+(0.018-0.012\mathrm i)z^2
 +\frac{0.035}{z-\zeta_1}-\frac{0.025\mathrm i}{z-\zeta_2},\label{eq:starfish-F}\\
 W(z)={}&-0.03\mathrm i+(0.045-0.020\mathrm i)z
 +(0.014+0.010\mathrm i)z^2
 +\frac{0.030\mathrm i}{z-\zeta_3}-\frac{0.022}{z-\zeta_4}.\label{eq:starfish-G}
\end{align}
The poles are
\begin{equation*}
 \zeta_1=2.50+1.00\mathrm i,\quad \zeta_2=-2.40+1.20\mathrm i,
 \qquad \zeta_3=2.30-1.40\mathrm i,\quad \zeta_4=-2.60-0.90\mathrm i.
\end{equation*}
All poles lie outside the starfish, so $F$ and $W$ are analytic in the interior.
With the stress $\Sigma$ of Section~\ref{sec:stokernels}, computed exactly from $F$ and $W$,
we prescribe on $\Gamma_S$
\begin{equation}\label{eq:starfish-boundary-data}
 q=u^-|_{\Gamma_S},\qquad f=\Sigma^-n.
\end{equation}
By the jump relations, the resulting velocity and pressure potentials equal
$u^-$ and $p^-$ inside $\Gamma_S$, vanish outside, and equal $u^-/2$ and
$p^-/2$ on $\Gamma_S$. The pressure densities $\phi$ and $\psi$ are computed by Fourier
differentiation, as for general data.

\begin{figure}[!htb]\centering
 \begin{minipage}[t]{.49\linewidth}\centering
  \includegraphics[width=\linewidth]{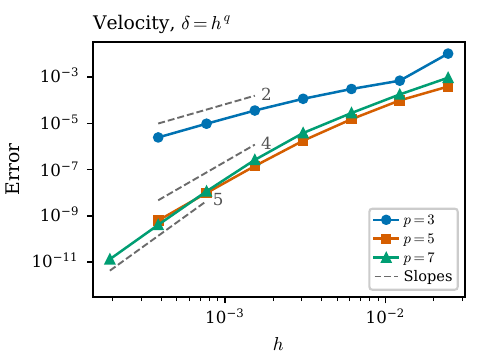}\par\vspace{-3pt}
  {\small (a)}
  \end{minipage}\hfill\begin{minipage}[t]{.49\linewidth}\centering
  \includegraphics[width=\linewidth]{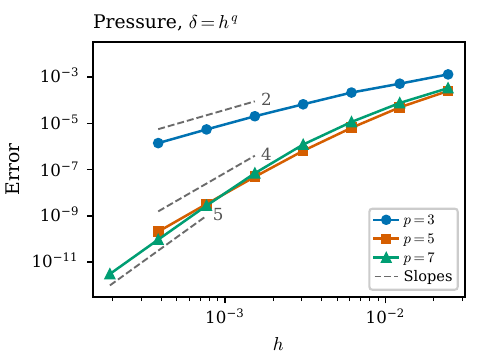}\par\vspace{-3pt}
  {\small (b)}
  \end{minipage}\par\smallskip
 \begin{minipage}[t]{.49\linewidth}\centering
  \includegraphics[width=\linewidth]{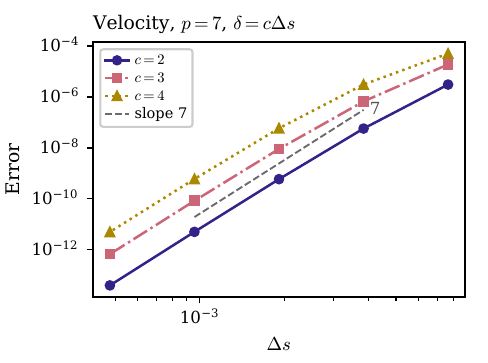}\par\vspace{-3pt}
  {\small (c)}
  \end{minipage}\hfill\begin{minipage}[t]{.49\linewidth}\centering
  \includegraphics[width=\linewidth]{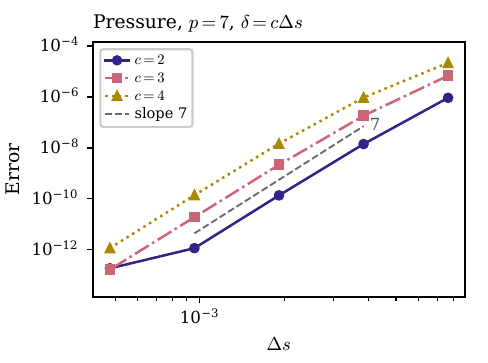}\par\vspace{-3pt}
  {\small (d)}
  \end{minipage}
 \caption{Stokes solution for the starfish: velocity (a),(c) and pressure
 (b),(d), with $\delta=h^q$ in (a),(b) and $p=7$, $\delta=c\Delta s$ in (c),(d).}
 \label{fig:stokes-starfish-results}
\end{figure}

We measure the maximum Euclidean norm of the velocity error and the
maximum pressure error. For $\delta$ as in
\eqref{eq:numscaling}, the velocity and pressure converge at close to the
predicted orders (Table~\ref{tab:rates},
Figure~\ref{fig:stokes-starfish-results}(a),(b)).

With the seventh-order kernel and $\delta=c \Delta s$,
$c=2,3,4$, the velocity error decreases at a rate close to seven
(Table~\ref{tab:rates}, Figure~\ref{fig:stokes-starfish-results}(c)). The
pressure error decreases at a similar rate over $N_\Gamma=2^{11}$--$2^{13}$
(Table~\ref{tab:rates}, Figure~\ref{fig:stokes-starfish-results}(d)); over
this range the velocity rates are $6.75$, $6.49$, and $6.17$. At
$N_\Gamma=2^{14}$ the pressure error levels off: it is limited by roundoff and
by the Fourier differentiation of the pressure densities
(Section~\ref{sec:pressure-densities}), and its estimated discretization error
is no longer well below the total error.

\subsubsection{The Stokeslet and stresslet integrals and the local corrections}\label{sec:correction-test}
To verify the accuracy of each layer without possible error cancellation between them,
we test the Stokeslet and stresslet integrals separately on the starfish. Exact integral
values are unavailable, so we compare with numerical reference values, prescribing smooth
densities independently as
\[
 \begin{aligned}
 q(t)&=\bigl(\sin(2t)+0.25\cos(5t),\;0.6\cos(3t)-0.2\sin t\bigr),\\
 f(t)&=\bigl(0.7\cos(2t)+0.15\sin(4t),\;0.5\sin(3t)+0.2\cos t\bigr).
 \end{aligned}
\]
The targets have the 49 closest-point parameters $t_0=2\pi(j+0.37)/49$, $j=0,\dots,48$,
and $b=\pm\delta/2$. We halve $\delta$ from $0.04$ to $0.00125$, with
$N_\Gamma=\max\{2^{13},2^{\lceil\log_2(32/\delta)\rceil}\}$. Reference values are computed
from the unregularized kernels by Gauss--Legendre quadrature on panels graded by bisection,
with singularity subtraction for the stresslet; the 48-point and 64-point rules agree to
within $10^{-14}$.

\begin{figure}[!htbp]
\centering
\begin{minipage}[t]{.49\linewidth}\centering
  \includegraphics[width=\linewidth]{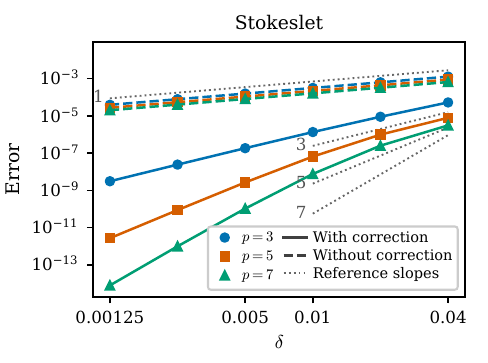}\par\vspace{-3pt}
  {\small (a)}
  \end{minipage}\hfill\begin{minipage}[t]{.49\linewidth}\centering
  \includegraphics[width=\linewidth]{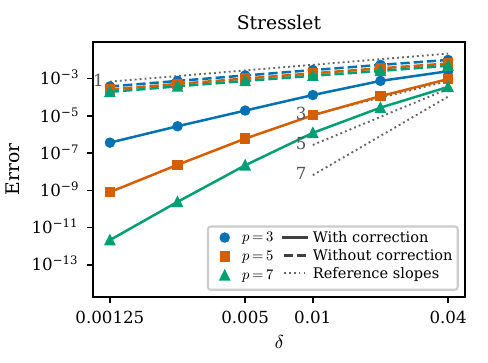}\par\vspace{-3pt}
  {\small (b)}
  \end{minipage}
\caption{Stokeslet and stresslet integrals for the starfish at
$b=-\delta/2$, with (solid) and without (dashed) the local corrections.}\label{fig:correction-ablation}
\end{figure}

Without the local corrections, the errors of both integrals are only about first order in
$\delta$ (Figure~\ref{fig:correction-ablation}, dashed curves). With them, the rates between
$\delta=0.0025$ and $0.00125$ are close to $p$: between $2.92$ and $3.04$ for $p=3$, $4.85$ and $5.05$
for $p=5$, and $6.72$ and $7.06$ for $p=7$, for both integrals and both sides of the curve.
For the seventh-order Stokeslet integral, whose error at $\delta=0.00125$ is at the level of
the reference error, the rate is taken between $\delta=0.005$ and $0.0025$.

\subsubsection{Differentiation of the pressure densities}\label{sec:pressure-densities}
For the exact starfish solution of Section~\ref{sec:num-starfish}, with the
boundary data \eqref{eq:starfish-boundary-data}, the pressure
densities on $\Gamma_S$ are also given by
\[
 \phi=-\partial_n p^-,\qquad \psi=-p^-.
\]
Indeed, incompressibility gives $\psi=-p^-$, while the Stokes equations
give $\phi=\mu\partial_s\omega^-=-\partial_n p^-$, with
$\omega^-=\partial_{x_1}u_2^--\partial_{x_2}u_1^-$ and the outward normal
convention used here. At $N_\Gamma=2^{10}$, the densities computed by
Fourier differentiation satisfy the identities for $\psi$ and $\phi$ to
within $10^{-13}$ and $10^{-10}$, respectively; the residuals grow with $N_\Gamma$ as roundoff in the
differentiation is amplified. To see how this affects the pressure, we
compute it once with the Fourier densities and once with the exact ones,
with everything else the same.

We use the 33 closest-point parameters $2\pi j/33$, each also shifted by half a
mesh interval, $b/\delta=-1/2,0,1/2$, the seventh-order kernel, and
$\delta=c\Delta s$ with $c=2,4$. (These targets differ from those of
Section~\ref{sec:num-starfish}.) At $N_\Gamma=2^{14}$ and $\delta=2\Delta s$, the
exact densities reduce the maximum pressure error by a factor of about $15$,
from about $10^{-13}$; on coarser meshes or with $c=4$, the regularization error
dominates and the two computations give nearly the same error. Thus the smallest
pressure errors are limited by the numerical differentiation of the densities.

\subsection{Targets approaching the curve}\label{sec:num-neartrace}
\begin{figure}[!htbp]
\centering
\includegraphics[width=.5\linewidth]{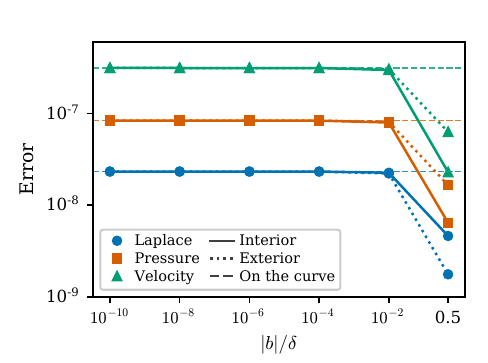}
\caption{Maximum errors as targets approach the curve: Laplace solution for the
kite, and Stokes velocity and pressure for the starfish.}
\label{fig:near-curve}
\end{figure}

The estimates of Sections~\ref{sec:lapbound} and~\ref{sec:stokesbounds} are
uniform for $|b|\le\Lambda\delta$; here we check that the computed values tend
to their limits on the curve as $b\to0$ from either side. We use the 33
closest-point parameters $t_0=2\pi(j+0.37)/33$, $j=0,\dots,32$,
$\delta=0.01$, and $N_\Gamma=2^{13}$. With the seventh-order kernels and
$|b|/\delta=10^{-10}$, the Stokeslet and stresslet integrals of
Section~\ref{sec:correction-test} differ from their values at $b=0$ by less
than $10^{-11}$; for the
stresslet, the value at $b=0$ is computed with $\chi=1$ or $0$, according to
the side of the target.
For the exact solutions on the kite and the starfish, with $p=7$ and
$|b|/\delta=0.5$, $10^{-2}$, $10^{-4}$, $10^{-6}$, $10^{-8}$, and $10^{-10}$
on both sides, the errors approach those at $b=0$ as $|b|$ decreases
(Figure~\ref{fig:near-curve}).

\subsection{Varying \texorpdfstring{$\delta$}{delta} with fixed nodes and targets}\label{sec:width-selection-new}
\begin{figure}[!htbp]
\centering
\begin{minipage}[t]{.49\linewidth}\centering
  \includegraphics[width=\linewidth]{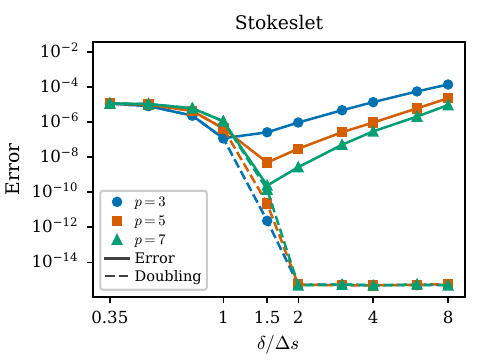}\par\vspace{-3pt}
  {\small (a)}
  \end{minipage}\hfill\begin{minipage}[t]{.49\linewidth}\centering
  \includegraphics[width=\linewidth]{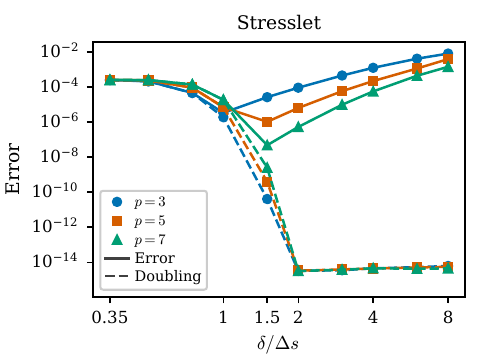}\par\vspace{-3pt}
  {\small (b)}
  \end{minipage}
\caption{Stokeslet and stresslet integrals for the starfish against
$\delta/\Delta s$: errors (solid) and estimated discretization errors (dashed).}
\label{fig:width-selection-new}
\end{figure}

To study the choice of $\delta$ separately from the target locations, we fix
$N_\Gamma=2^{11}$ and $b=-0.003$ for the Stokeslet and stresslet integrals of
Section~\ref{sec:correction-test} and vary $\delta/\Delta s$ from $0.35$ to $8$
(Figure~\ref{fig:width-selection-new}). For small $\delta$ the estimated discretization
error is comparable with the total error; for large $\delta$ it is negligible, but the
regularization error is larger. For these targets the smallest fifth- and seventh-order
errors occur near $\delta/\Delta s=1.5$, but there the estimated discretization error of
the seventh-order Stokeslet integral is comparable with the total error. At $\delta/\Delta s=2$, the estimated discretization errors of both
integrals are at roundoff level.
In practice, one can decrease $\delta$ while monitoring the estimated discretization error;
once it becomes significant, decreasing $\delta$ further requires more nodes. The product
$\delta\|\kappa\|_\infty$ must also be small, as Proposition~\ref{prop:lap} assumes. The best
ratio $\delta/\Delta s$ depends on the geometry, the density, the order, and the target
distance.

\FloatBarrier

\section{Values on a Cartesian grid}\label{sec:num-allgrid}
In many applications the solution is needed on a whole grid rather than at a
few targets, for example when the flow is coupled with an equation for a heat
or solute concentration in the surrounding fluid, or when the layer potentials
supply the interface data for a grid-based method. Evaluating the integrals at
every grid point would be costly. Following Mayo \cite{Mayo1984,Mayo1985} and
\cite{BL2001,BT2026}, we evaluate the regularized potentials directly only at
grid points near $\Gamma$ and on the outer boundary, and obtain the values
elsewhere with a fast Poisson solver based on the discrete sine transform (DST).

\subsection{Finite-difference discretization}\label{sec:grid-disc}
Let $h_c$ be the grid spacing of a square computational box and $M$ the
number of interior grid points in each direction.
We discretize the Laplacian by the nine-point compact operator
\cite{Iserles2008}
\begin{equation}\label{eq:mehrstellen}
 \begin{aligned}
 (L_{h_c}v)_{ij}=\frac{1}{6h_c^2}\big[&4(v_{i+1,j}+v_{i-1,j}
 +v_{i,j+1}+v_{i,j-1})\\[-2pt]
 &+v_{i+1,j+1}+v_{i+1,j-1}+v_{i-1,j+1}+v_{i-1,j-1}-20v_{ij}\big].
 \end{aligned}
\end{equation}
For smooth $v$, Taylor expansion gives
\begin{equation}\label{eq:mehrstellenexp}
 L_{h_c}v=\Delta v+\frac{h_c^2}{12}\Delta^2v
 +\frac{h_c^4}{360}\big(
 \partial_{x_1}^6v+5\partial_{x_1}^4\partial_{x_2}^2v
 +5\partial_{x_1}^2\partial_{x_2}^4v+\partial_{x_2}^6v\big)+O(h_c^6).
\end{equation}
The operator is second order accurate for a general smooth function, and
its truncation error is $O(h_c^6)$ for a harmonic function. For the Stokes
velocity,
$\Delta u=\mu^{-1}\nabla p$ and $\Delta p=0$, hence $\Delta^2u=0$, and the
truncation error is $O(h_c^4)$.

We approximate the pressure gradient by the centered fourth-order formulas
\begin{equation}\label{eq:grad4}
 \begin{aligned}
 (D_{x_1}^{(4)}v)_{ij}&=
 \frac{v_{i-2,j}-8v_{i-1,j}+8v_{i+1,j}-v_{i+2,j}}{12h_c},\\
 (D_{x_2}^{(4)}v)_{ij}&=
 \frac{v_{i,j-2}-8v_{i,j-1}+8v_{i,j+1}-v_{i,j+2}}{12h_c}.
 \end{aligned}
\end{equation}
We write $\nabla_{h_c}^{(4)}=(D_{x_1}^{(4)},D_{x_2}^{(4)})$, with fourth-order
one-sided formulas next to the outer boundary.

\subsection{Irregular points and directly evaluated values}\label{sec:grid-irregular}
\begin{figure}[tb]\centering
 \includegraphics[width=.64\linewidth]{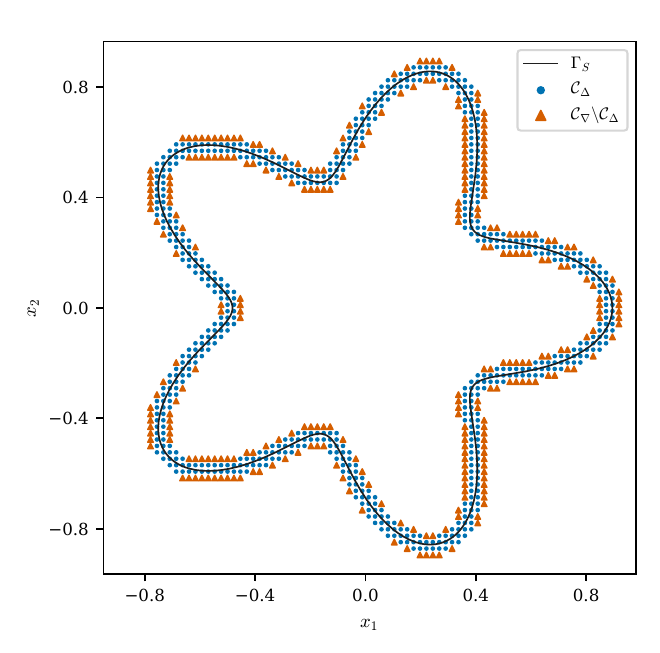}
 \caption{Irregular points for the starfish at $M=2^{7}$: the set
 $\mathcal{C}_\Delta$ and the additional points of
 $\mathcal{C}_\nabla\setminus\mathcal{C}_\Delta$, where only the
 pressure-gradient stencil crosses $\Gamma$.}
 \label{fig:allgrid-stencils}
\end{figure}

Let $\mathcal{W}_\Delta(x_{ij})$ consist of $x_{ij}$, its four axial neighbors,
and its four diagonal neighbors, as in \eqref{eq:mehrstellen}, and let
$\mathcal{W}_\nabla(x_{ij})$ consist of $x_{ij}$ and the points at axial offsets
$\pm h_c,\pm2h_c$, as in \eqref{eq:grad4}. A stencil crosses $\Gamma$ if its
points do not all lie on the same side of $\Gamma$; whether a grid point is
inside or outside $\Omega$ is found by a point-in-polygon test. The sets of irregular
points \cite{BLay2006,LL1994} are
\begin{equation}\label{eq:crossingsets}
 \mathcal{C}_\Delta=\{x_{ij}:\mathcal{W}_\Delta(x_{ij})\text{ crosses }\Gamma\},
 \qquad
 \mathcal{C}_\nabla=\{x_{ij}:\mathcal{W}_\nabla(x_{ij})\text{ crosses }\Gamma\}.
\end{equation}
Grid points outside $\mathcal{C}_\Delta$ are regular for the Laplace and
pressure problems. For the Stokes velocity, the irregular points are
\begin{equation}\label{eq:velocitycrossing}
 \mathcal{C}_u=\mathcal{C}_\Delta\cup\mathcal{C}_\nabla,
\end{equation}
since the pressure-gradient stencil can cross $\Gamma$ when the nine-point
stencil does not (Figure~\ref{fig:allgrid-stencils}).

The layer potentials are evaluated directly at the irregular points and at
all points of their nine-point stencils, that is, on
\begin{equation}\label{eq:halosets}
 \mathcal{E}_L=\bigcup_{x\in\mathcal{C}_\Delta}\mathcal{W}_\Delta(x),\qquad
 \mathcal{E}_p=\mathcal{E}_u=\bigcup_{x\in\mathcal{C}_u}\mathcal{W}_\Delta(x).
\end{equation}
The pressure equation needs directly evaluated values only on
$\mathcal{E}_L$, but we evaluate the pressure directly, and keep these values,
on the larger set $\mathcal{E}_p=\mathcal{E}_u$, where the velocity is also
evaluated directly.

\subsection{Laplace problem}
For a piecewise harmonic function $v=Sf+Dg$, we evaluate
$v^I=S_\delta^{(p)}f+D_\delta^{(p)}g$ on $\mathcal{E}_L$, taking the interior or
exterior value according to the side of each grid point. We set
\begin{equation}\label{eq:lapgridrhs}
 F_h^L(x_{ij})=\begin{cases}
 L_{h_c}v^I(x_{ij}),&x_{ij}\in\mathcal{C}_\Delta,\\
 0,&x_{ij}\notin\mathcal{C}_\Delta,
 \end{cases}
\end{equation}
and solve $L_{h_c}v_h=F_h^L$ with Dirichlet data on the outer boundary
computed from the layer potentials with the original kernels, which are not
singular there. The reported solution uses the directly evaluated values
on $\mathcal{E}_L$ and the grid solution elsewhere:
\begin{equation}\label{eq:laphybrid}
 v_h^{\rm all}(x)=\begin{cases}v^I(x),&x\in\mathcal{E}_L,\\
 v_h(x),&x\notin\mathcal{E}_L.
 \end{cases}
\end{equation}

To relate the error of the grid solution to that of the directly evaluated
values, let $e^I=v^I-v$ and $e_h=v_h-v$, and let $\tau_h=L_{h_c}v$ be the
truncation error at the regular points, where $\Delta v=0$. With exact
boundary data,
\begin{equation}\label{eq:reconstruction-error-equation}
 L_{h_c}e_h
 =\mathbf1_{\mathcal{C}_\Delta}L_{h_c}e^I
  -\mathbf1_{\mathcal{C}_\Delta^c}\tau_h.
\end{equation}
Thus $L_{h_c}$ is applied to the error of the directly evaluated values
before the solve, and a bound on this error in the maximum norm alone
allows growth by a factor $h_c^{-2}$ in the right-hand side. This term,
however, is supported on $\mathcal{C}_\Delta$, within $O(h_c)$ of $\Gamma$.
For the standard discrete Laplacian, Beale and Layton
\cite[Lemmas~2.2 and~2.3]{BLay2006} show that such a grid function is the
discrete divergence of a grid function smaller by a factor $O(h_c)$, and
that the discrete solution is bounded by the maximum norm of this function
plus that of the truncation error, uniformly in $h_c$; its first
differences are bounded by the same quantity times $|\log h_c|$. If these
estimates also hold for $L_{h_c}$, as Beale and Tlupova expect for their
fourth-order Laplacian \cite[Sect.~4.3]{BT2026}, then
\begin{equation}\label{eq:grid-heuristic}
 \max|e_h|\le C\Bigl(h_c^{-1}\max_{\mathcal{E}_L}|e^I|
  +\max_{\mathcal{C}_\Delta^c}|\tau_h|\Bigr),
\end{equation}
and an error $O(h_c^5)$ in the directly evaluated values gives an error
$O(h_c^4)$ in the grid solution. For Stokes, the pressure error is also
differentiated before the velocity solve. The fourth-order gradient is a
combination of first differences, so the same argument gives a velocity
error $O(h_c^4|\log h_c|)$.

\subsection{Stokes problem: pressure first, then velocity}
We first evaluate
$p^I=S_\delta^{(p)}\phi-D_\delta^{(p)}\psi$ on $\mathcal{E}_p$ and solve
\begin{equation}\label{eq:pressuregridrhs}
 L_{h_c}p_h=F_h^p,\qquad
 F_h^p(x_{ij})=\begin{cases}
 L_{h_c}p^I(x_{ij}),&x_{ij}\in\mathcal{C}_\Delta,\\
 0,&x_{ij}\notin\mathcal{C}_\Delta.
 \end{cases}
\end{equation}
As in \eqref{eq:laphybrid}, the reported pressure $p_h^{\rm all}$ is $p^I$
on $\mathcal{E}_p$ and $p_h$ elsewhere. At the regular points for the
velocity, we differentiate $p_h$, not $p_h^{\rm all}$.

We then evaluate the corrected Stokes velocity $u^I$ on $\mathcal{E}_u$ and set
\begin{equation}\label{eq:velocitygridrhs}
 F_h^u(x_{ij})=\begin{cases}
 L_{h_c}u^I(x_{ij}),&x_{ij}\in\mathcal{C}_u,\\[2pt]
 \mu^{-1}\nabla_{h_c}^{(4)}p_h(x_{ij}),
 &x_{ij}\notin\mathcal{C}_u.
 \end{cases}
\end{equation}
We solve $L_{h_c}(u_h)_k=(F_h^u)_k$ for the two velocity components, with
Dirichlet data on the outer boundary. The reported velocity $u_h^{\rm all}$
is $u^I$ on $\mathcal{E}_u$ and $u_h$ elsewhere.

The direct evaluations cost $O(h_c^{-2})$, since there are
$N_\Gamma=O(h_c^{-1})$ nodes on the curve and $O(h_c^{-1})$ grid points where
the layer potentials are evaluated directly, and each DST solve costs
$O(h_c^{-2}\log(1/h_c))$, whereas direct
evaluation at all $O(h_c^{-2})$ grid points would cost $O(h_c^{-3})$. The
direct evaluations could be accelerated by the fast multipole method
\cite{GR1987} or a treecode \cite{WKT2020}; we do not use them here.

\subsection{Results on a Cartesian grid}\label{sec:grid-results}
We apply this method to the Laplace solution
for the kite of Section~\ref{sec:num-kite} and to the Stokes solution for the
starfish of Section~\ref{sec:num-starfish}, on the box $(-1.5,1.5)^2$, so
that $h_c=3/(M+1)$. We use the seventh-order kernel with
$\delta=h^{5/7}$, as in \eqref{eq:numscaling}. For each curve and
each $M=2^{6},2^{7},\dots,2^{10}$, $N_\Gamma$ is the smallest even integer
at least $2\pi\max_t|x'(t)|/h_c$ with no prime factor larger than $11$, for
the FFT. The ratio $\Delta s/h_c$, with
$\Delta s$ as in \eqref{eq:numphysical}, then lies
between $0.96$ and $1.00$.

Since $h\asymp \Delta s\approx h_c$, the regularization error at the
directly evaluated points is $O(\delta^7)=O(h_c^5)$. The truncation errors
(Section~\ref{sec:grid-disc}) are $O(h_c^6)$ for the harmonic Laplace and
pressure fields and $O(h_c^4)$ for the velocity. If \eqref{eq:grid-heuristic}
holds, the grid solution therefore has error $O(h_c^4)$, or
$O(h_c^4|\log h_c|)$ for the velocity.
In addition to the maximum error, we report the discrete $L^2$ norm
\begin{equation}\label{eq:discreteL2}
 \|e\|_{2,h}=\left(h_c^2\sum_{i,j=0}^{M+1}|e_{ij}|^2\right)^{1/2}.
\end{equation}
Here the sum includes the grid points on the outer boundary and, for the
velocity, $|e_{ij}|$ is the Euclidean norm of the error.

\begin{figure}[!htb]\centering
 \begin{minipage}[t]{.49\linewidth}\centering
  \includegraphics[width=\linewidth]{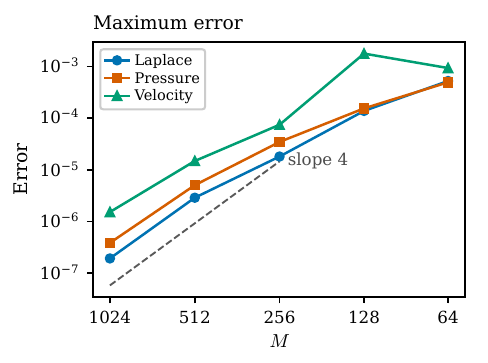}\par\vspace{-3pt}
  {\small (a)}
  \end{minipage}\hfill
 \begin{minipage}[t]{.49\linewidth}\centering
  \includegraphics[width=\linewidth]{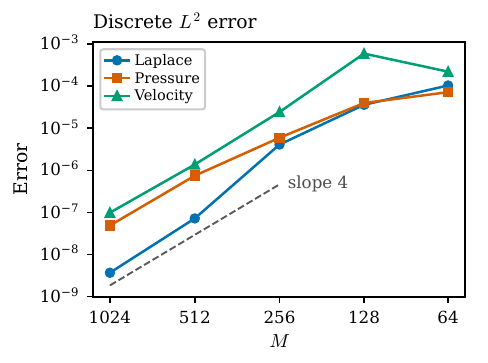}\par\vspace{-3pt}
  {\small (b)}
  \end{minipage}
 \caption{Maximum (a) and discrete $L^2$ (b) errors on the Cartesian grid for the
 Laplace solution for the kite and the Stokes solution for the starfish.}
 \label{fig:allgrid-results}
\end{figure}

Figure~\ref{fig:allgrid-results}(a),(b) shows the maximum and discrete
$L^2$ errors. On every grid, all three maximum errors occur at points where the potentials
are evaluated directly, so the maximum norm mainly measures the error of
these values near the interface, while the discrete $L^2$ norm also reflects
the grid solution.

\begin{figure}[!htb]\centering
 \includegraphics[width=.32\linewidth]{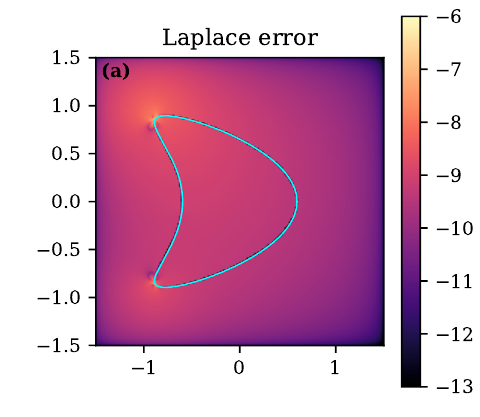}\hfill
 \includegraphics[width=.32\linewidth]{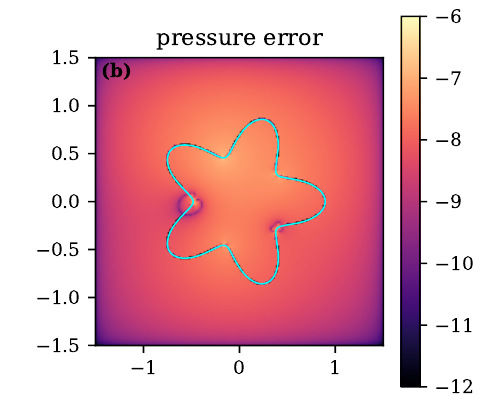}\hfill
 \includegraphics[width=.32\linewidth]{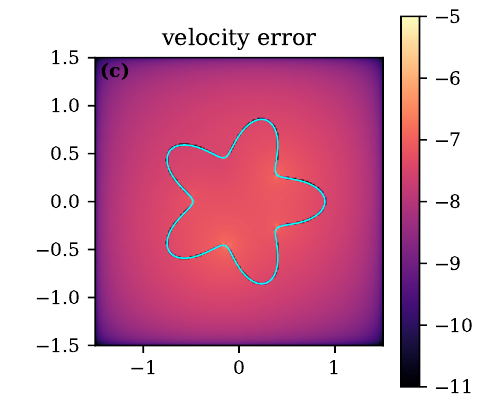}
 \caption{$\log_{10}$ of the error on the Cartesian grid at $M=2^{10}$: (a) Laplace
 solution for the kite; (b) pressure and (c) velocity for the starfish.}
 \label{fig:allgrid-maps}
\end{figure}

Between $M=2^{9}$ and $2^{10}$, the rates range from $3.91$ for the Laplace
solution to $3.29$ for the velocity in the maximum norm, and from $4.30$ to
$3.81$ in the discrete $L^2$ norm.
The velocity error at $M=2^{7}$ is anomalously large; we return to this in
Section~\ref{sec:grid-curvature}.
On the finest grid, the largest errors occur near the points of highest
curvature, the two tips of the kite and the inner corners between the arms
of the starfish (Figure~\ref{fig:allgrid-maps}).

\subsection{Exact values at the directly evaluated points}\label{sec:grid-exact}
To isolate the error of the grid method, we repeat the computation with
exact values at the directly evaluated points and on the outer boundary,
leaving everything else unchanged; for Stokes, the pressure is still
computed first and differentiated with $\nabla_{h_c}^{(4)}$. At $M=2^{9}$ the
maximum errors are below $10^{-11}$, many orders of magnitude smaller
than those of the full method. With the location of the maximum errors
noted in Section~\ref{sec:grid-results}, this indicates that the regularized
evaluations near the interface determine the maximum error.

To measure the error of the grid method above roundoff, we repeat this
test with exact solutions that vary on a shorter length scale. Let
\[
 \mathcal{A}(x_1,x_2)=e^{k(x_1-0.895)}\cos(kx_2),\qquad
 \mathcal{B}(x_1,x_2)=e^{k(x_1-0.895)}\sin(kx_2),\qquad k=8.
\]
Since $\max_t\|x(t)\|_\infty=0.895$ for both curves, the shift by $0.895$
gives $|\mathcal{A}|,|\mathcal{B}|\le1$ inside them. For Laplace, we take
$v^-=\mathcal{A}$ inside the kite and $v^+=0$ outside.
For Stokes, with $\mu=1$, the stream function $\Psi=x_2\mathcal{A}/k$
inside the starfish gives
\[
 u^-=(\partial_{x_2}\Psi,-\partial_{x_1}\Psi)
      =(\mathcal{A}/k-x_2\mathcal{B},-x_2\mathcal{A}),\qquad p^-=-2\mathcal{A},
\]
and we set the exterior velocity and pressure to zero. These fields satisfy
$\nabla\cdot u^-=0$, $\Delta u^-=\nabla p^-$, and $\Delta p^-=0$.

\begin{figure}[htbp]\centering
 \begin{minipage}[t]{.49\linewidth}\centering
  \includegraphics[width=\linewidth]{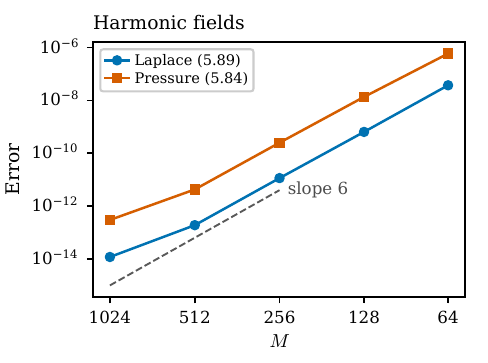}\par\vspace{-3pt}
  {\small (a)}
  \end{minipage}\hfill\begin{minipage}[t]{.49\linewidth}\centering
  \includegraphics[width=\linewidth]{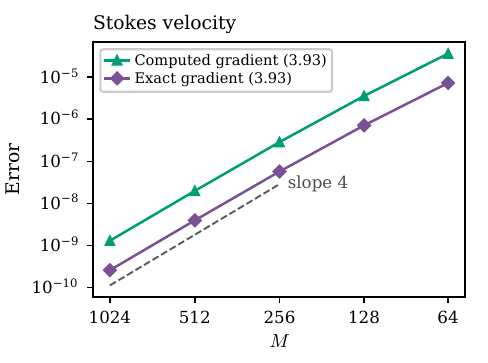}\par\vspace{-3pt}
  {\small (b)}
  \end{minipage}
 \caption{Cartesian grid method with exact values at the directly
 evaluated points: (a) Laplace solution and pressure; (b) velocity with computed
 and exact pressure gradients.}
 \label{fig:exact-data-reconstruction}
\end{figure}

The errors decrease at rates close to six for the harmonic Laplace and
pressure fields over $M=2^{7}$--$2^{9}$, the finer grids being affected by
roundoff, and close to four for the velocity between $M=2^{9}$ and $2^{10}$, with
both the computed and the exact pressure gradient
(Figure~\ref{fig:exact-data-reconstruction}); the exact
gradient reduces the velocity error by a factor of about five. These rates
agree with the truncation errors of \eqref{eq:mehrstellenexp} for harmonic
and biharmonic functions.

\subsection{The choice of \texorpdfstring{$\delta$}{delta} and the curvature}\label{sec:grid-curvature}
\begin{figure}[htbp]
\centering
\includegraphics[width=.5\linewidth]{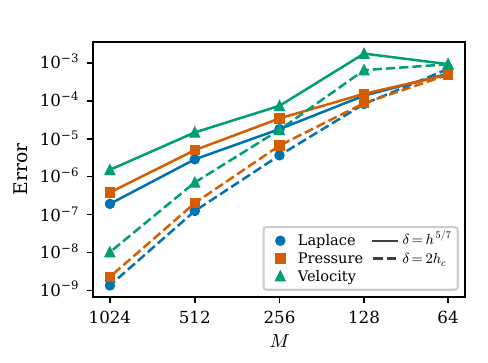}
\caption{Maximum errors on the Cartesian grid for $\delta=h^{5/7}$ (solid) and
$\delta=2h_c$ (dashed).}
\label{fig:allgrid-widths}
\end{figure}

We repeat the computations with $\delta=2h_c$, proportional to the mesh
size, on the same grids and with the same nodes on the curves. On the finer
grids, $\delta=2h_c$ is much more accurate than $\delta=h^{5/7}$
(Figure~\ref{fig:allgrid-widths}). At $M=2^{6}$ the two values of $\delta$
differ by less than 7\% and give similar errors. At
$M=2^{10}$, $\delta=2h_c$ reduces the errors by about two orders of
magnitude. At about 100 targets on this grid, including the points with the
largest errors, the estimated discretization errors are below $10^{-13}$.
Because the estimate of
Section~\ref{sec:discrete-error} does not imply convergence when $\delta$ is
proportional to the mesh size, this comparison does not indicate the
asymptotic rate.

On coarse grids, the velocity error is also affected by the geometric
factors in the local correction. Table~\ref{tab:geometry-resolution} gives,
for the starfish, $\min_{y\in\mathcal{E}_u}(1+b\kappa_0)$, which does not
depend on $\delta$, and $\delta\|\kappa\|_\infty$, with
$\|\kappa\|_\infty\approx20.2$.
\begin{table}[htbp]\centering
\caption{Geometric quantities for the Stokes computation for the starfish.}
\label{tab:geometry-resolution}
\begin{tabular}{rccc}
\toprule
$M$ & $\min_{\mathcal{E}_u}(1+b\kappa_0)$ &
$h^{5/7}\|\kappa\|_\infty$ & $2h_c\|\kappa\|_\infty$\\
\midrule
$2^{6}$   & $0.32$ & $1.86$ & $1.86$\\
$2^{7}$  & $0.21$ & $1.14$ & $0.94$\\
$2^{8}$  & $0.44$ & $0.71$ & $0.47$\\
$2^{9}$  & $0.65$ & $0.44$ & $0.24$\\
$2^{10}$ & $0.83$ & $0.27$ & $0.12$\\
\bottomrule
\end{tabular}
\end{table}
At $M=2^{7}$, the largest velocity error occurs at a point with
$\alpha^2=1+b\kappa_0\approx0.21$. There the factors $\alpha^{-3}$ and
$\alpha^{-5}$ in \eqref{eq:C2} make the stresslet correction large, while the
discretization error and the error in the closest point are negligible. With
$\delta=h^{5/7}$ on this grid, $\delta\|\kappa\|_\infty\approx1.14$, so $\delta$
exceeds the smallest radius of curvature, contrary to the hypothesis
$\Lambda\delta\|\kappa\|_\infty\le\frac12$ of Proposition~\ref{prop:stokes}. On
the finest grid $\delta\|\kappa\|_\infty\approx0.27$ and $\alpha$ is close to one
(Table~\ref{tab:geometry-resolution}).
\FloatBarrier

\section{Conclusions}\label{sec:conclusions}
We have constructed regularized versions of the kernels for Laplace and
Stokes layer potentials in two dimensions. They allow control of the
discretization error: we prove that the regularization error is
$O(\delta^p)$ in the smoothing radius $\delta$, $p=3,5,7$, uniformly for
targets on and near the curve, so $\delta$ can be taken several times the
node spacing $h$ to keep the trapezoidal discretization error small. The
method is easy to implement and use: for every target, the trapezoidal rule
is applied on one standard set of nodes on the curve, reusing the same
samples of the curve and density.

The kernel coefficients are explicit functions of the scaled distance
$b/\delta$ to the curve, and the Laplace potentials require no derivatives
of the density. The Stokes kernels use the same coefficients together with
explicit local corrections that remove one remaining error moment, involve
the curvature, and vanish on the curve; without them, the Stokeslet and
stresslet errors are only first order in $\delta$. The stresslet correction
and the pressure formula use tangential derivatives of the densities.

The error bounds hold when $\delta$ is at most a fixed fraction of the
smallest radius of curvature. For $\delta=ch^q$ with $q<1$, the
discretization error decays faster than $\delta^p$, and the observed
convergence rates are close to $pq$. With $p=7$ and $\delta$ equal to two
to four times the largest distance between adjacent nodes, we observe rates
of $6.4$ to $7$ while the discretization error remains small; doubling the
number of nodes at fixed $\delta$ estimates the discretization error.

On a Cartesian grid, a fast Poisson solver supplies values at all grid
points from the ones near the curve whose values are computed directly,
where the largest errors occur. On finer grids, taking $\delta$ equal to
twice the grid spacing is much more accurate than $\delta=h^{5/7}$ once
$\delta$ resolves the curvature. Future work includes choosing $\delta$
locally along the curve from the node spacing or curvature, solving
boundary integral equations, and treating open arcs and corners.

\appendix
\section{Evaluation of the kernels and coefficients}\label{app:implementation}
The formulas below avoid loss of accuracy from cancellation near $\xi=0$ and for large $|\lambda|$.

\subsection{The logarithmic terms}
When the source and target are close, we evaluate $2\log r+E_1(\xi)$ as a
single expression, since each term diverges as $r\to0$ while their sum is
analytic. Substituting the expansion \cite[Chap.~6]{NIST2010}
\begin{equation}\label{eq:E1series}
 E_1(\xi)=-\gamma-\log\xi+\sum_{k=1}^{\infty}\frac{(-1)^{k+1}\xi^{k}}{k\,k!},
\end{equation}
we obtain
$2\log r+E_1(\xi)=2\log\delta-\gamma+\xi-\xi^2/4+\xi^3/18-\cdots$, which
contains no singular terms. For $\xi<1/2$, we evaluate the sum from the
first 18 terms of this series; otherwise, we evaluate $E_1$ with the
function \texttt{exp1} of SciPy \cite{SciPy2020}; the same function is available
as \texttt{expint} in MATLAB and as \texttt{gsl\_sf\_expint\_E1} in the GNU
Scientific Library.

\subsection{Coefficients}\label{app:coeff}
The functions $\mathcal{H}_n$ satisfy a two-term recurrence. With
$\erfcx(x)=e^{x^2}\operatorname{erfc}(x)$,
\begin{equation}\label{eq:Hrec}
 \mathcal{H}_0=2\bigl[1-\sqrt\pi\,|\lambda|\,\erfcx|\lambda|\bigr],\qquad
 \mathcal{H}_n=\frac{2}{n+1}\bigl(1-\lambda^2\mathcal{H}_{n-2}\bigr),
 \quad n=2,4,6,\ldots .
\end{equation}
Thus all $\mathcal{H}_n$ are determined by $\mathcal{H}_0$.

The recurrence is accurate for moderate $|\lambda|$, but loses digits through
cancellation for large $|\lambda|$. The differences
$\mathcal{H}_0-\mathcal{H}_2$ and $\mathcal{H}_0-2\mathcal{H}_2+\mathcal{H}_4$ in
\eqref{eq:c5} and \eqref{eq:c7} decay faster than the individual $\mathcal{H}_n$
and are even more affected. We therefore compute these differences
directly from integrals with positive integrands, obtained from
\eqref{eq:Hdef}:
\begin{equation}\label{eq:Dint}
 \mathcal{U}_j(\lambda)=\lambda^{-2j-2}\int_0^\infty e^{-v}v^{j}
 \Bigl(1+\frac{v}{\lambda^2}\Bigr)^{-j-3/2}\dd v,
 \qquad j=0,1,2,\quad \lambda\ne0,
\end{equation}
for which $\mathcal{U}_0=\mathcal{H}_0$, $\mathcal{U}_1=\mathcal{H}_0-\mathcal{H}_2$, and
$\mathcal{U}_2=\mathcal{H}_0-2\mathcal{H}_2+\mathcal{H}_4$. The integrand is bounded,
and the small factor $\lambda^{-2j-2}$ is applied after the integration. We
use these formulas for $|\lambda|>1$ and the recurrence \eqref{eq:Hrec} for
$|\lambda|\le1$, with $\erfcx$ in $\mathcal{H}_0$ so that the product of
$e^{\lambda^2}$ and $\operatorname{erfc}|\lambda|$ is never formed. For
$|\lambda|>1$, the integrals in \eqref{eq:Dint} are computed with the
adaptive QUADPACK routine for infinite intervals \cite{QUADPACK1983}, called through
\texttt{scipy.integrate.\allowbreak quad} \cite{SciPy2020}, with absolute and relative
tolerances $2\times10^{-14}$ and at most 50 subintervals.

\begin{remark}\label{rem:nomoment}
The analysis requires only bounded $\lambda$, so the behavior of
$a_{k}^{(p)}(\lambda)$ as $|\lambda|\to\infty$ is not needed. For accuracy that is
uniform in the target distance, the coefficients are computed at the actual
scaled distance $\lambda=b/\delta$ of each target.
\end{remark}

\subsection{Radial factors in the Stokes kernels}
To avoid cancellation near $\xi=0$, we rewrite the radial factors
$A=A^{(p)}$ and $B=B^{(p)}=-A'$ of \eqref{eq:Udelta} and \eqref{eq:Hdelta},
for fixed $p$ and $\lambda$, with primes denoting $\partial_\xi$. Write
$P^{(p)}(\xi;\lambda)=1+\xi Z$, so that $Z$ is the polynomial
$\sum_{k=1}^m(-1)^ka_{k}^{(p)}(\lambda)\xi^{k-1}$ of degree $m-1$. Then
\begin{equation}\label{eq:radial}
 \begin{aligned}
 \Phi(\xi)&=\frac{1-e^{-\xi}}{\xi},\qquad A=\Phi-e^{-\xi}Z,\\
 \partial_\xi d^{(p)}&=e^{-\xi}\bigl(1+\xi Z-Z-\xi Z'\bigr),\qquad
 B=-A'=-\Phi'+e^{-\xi}(Z'-Z),
\end{aligned}
\end{equation}
so $A$, $B$, and $\partial_\xi d^{(p)}$ are computed from $\Phi$ and the
polynomial $Z$, without dividing a small difference by $\xi$. We compute
$\Phi$ as $-\operatorname{expm1}(-\xi)/\xi$, with $\Phi(0)=1$, where
$\operatorname{expm1}(z)=e^z-1$ is accurate for small $z$. For $\xi<1/2$,
including $\xi=0$, we compute $\Phi'$ from the first 20 terms of the series
$\Phi'(\xi)=\sum_{j\ge1}(-1)^jj\,\xi^{j-1}/(j+1)!$, which avoids the
cancellation between $\xi e^{-\xi}$ and $1-e^{-\xi}$; otherwise we use
$\Phi'(\xi)=[\xi e^{-\xi}+\operatorname{expm1}(-\xi)]/\xi^2$. The tensors are then
assembled as $r_ir_jA/\delta^2$ and
$r_ir_j(\boldsymbol r\cdot n)(-A')/\delta^4$, so neither $r^{-2}$ nor $r^{-4}$
is ever formed. These formulas give
$B(0)=\tfrac12+a_{1}^{(p)}+a_{2}^{(p)}$, in agreement with the limit in
Section~\ref{sec:stokernels}.

\section*{Acknowledgments}
The author thanks J.~Thomas Beale for explaining the problem and for many
helpful suggestions on the first draft of this article, and Yoichiro Mori for
hosting the visit of J.~Thomas Beale to the University of Pennsylvania during
which this work began.

\section*{Declarations}
\paragraph{Competing interests} The author declares no competing interests.

\paragraph{Use of AI tools} The author used Claude Code and Gemini to assist with
code generation and writing/editing the manuscript. After using these tools,
the author reviewed and verified all mathematical derivations, numerical code,
and text, and takes full responsibility for the content of the article.

\paragraph{Code and data availability} The code and the data for the
numerical examples are available in the GitHub repository
\url{https://github.com/han-zhou-math/high-order-reg-2d.git}.

\end{document}